\documentclass[ejsv2]{imsart}

\usepackage[authoryear,round,sort&compress]{natbib}
\setcitestyle{authoryear,round,semicolon}
\usepackage[T1]{fontenc}
\usepackage{newtxtext}

\usepackage[colorlinks,citecolor=blue,urlcolor=blue]{hyperref}

\DeclareMathOperator{\tr}{tr}
\usepackage{graphicx,subcaption}
\usepackage{amsmath,amsthm,mathtools,bm}
\usepackage{booktabs,multirow}
\usepackage{algorithm,algpseudocode}
\usepackage{tikz}

\newtheorem{theorem}    {Theorem}[section]
\newtheorem{lemma}      [theorem]{Lemma}
\newtheorem{proposition}[theorem]{Proposition}

\theoremstyle{definition}

\theoremstyle{remark}
\newtheorem{remark}[theorem]{Remark}
\theoremstyle{remark}
\newtheorem*{remark*}{Remark}  
\newcommand{\lat}{\mathrm{lat}}

\newtheorem{corollary}  [theorem]{Corollary} 
\newtheorem{assumption} [theorem]{Assumption}
\newcommand{\Z}{\mathbb Z}
\newcommand{\R}{\mathbb R}
\newcommand{\T}{\mathbb T}
\newcommand{\E}{\mathbb E}
\newcommand{\Var}{\operatorname{Var}}
\newcommand{\Cov}{\operatorname{Cov}}
\newcommand{\dd}{\mathrm d}

\newcommand{\Mat}{\mathrm{Mat}}

\usepackage{mathtools}
\DeclarePairedDelimiter{\abs}{\lvert}{\rvert}
\DeclarePairedDelimiter{\norm}{\lVert}{\rVert}

\newcommand{\PL}{\mathrm{PL}}   

\usetikzlibrary{calc,arrows.meta}
\tikzset{
  grid/.style={gray!30},
  stencil/.style={blue!70,thick,rounded corners=2pt},
  mini/.style={blue!60!black,dashed,rounded corners=1.5pt},
  coefplus/.style={fill=green!12,draw=green!50!black,rounded corners=1pt,inner sep=1.2pt},
  coefminus/.style={fill=red!10,draw=red!55!black,rounded corners=1pt,inner sep=1.2pt},
  sumop/.style={circle,draw=black,fill=gray!5,inner sep=1.2pt}
}

\begin{document}
\begin{frontmatter}

\title{Fixed and Increasing Domain Asymptotics for the Roughness and Scale of Isotropic Gaussian Random Fields}
\runtitle{Asymptotics for Moment Estimators}

\begin{aug}
\author{\fnms{Varun} \snm{Kotharkar}\ead[label=e1]{vsk34@stat.rutgers.edu}}
\and
\author{\fnms{Michael L.} \snm{Stein}\ead[label=e2]{ms2870@stat.rutgers.edu}}

\runauthor{V. Kotharkar and M. Stein}

\address{Department of Statistics, Rutgers University,
New Brunswick, NJ, USA. \printead{e1,e2}}
\end{aug}

\begin{abstract}
We establish a rigorous asymptotic theory for the joint estimation of
roughness and scale parameters in two-dimensional Gaussian random fields
with power-law generalized covariances \cite{Matheron1973, Stein1999, Yaglom1987}.
Our main results are bivariate central limit theorems for a class of method-of-moments estimators under increasing-domain and fixed-domain asymptotics.
The fixed-domain result follows immediately from the increasing-domain result from the self-similarity of Gaussian random fields with power-law generalized covariances \cite{IstasLang1997, Coeurjolly2001, ZhuStein2002}.
These results provide a unified distributional framework across these two classical regimes \cite{AvramLeonenkoSakhno2010-ESAIM, BiermeBonamiLeon2011-EJP} that makes the unusual behavior of the estimates under fixed-domain asymptotics intuitively obvious.
Our increasing-domain asymptotic results use spatial averages of quadratic forms of (iterated) bilinear product difference filters that yield explicit expressions for the estimates of roughness and scale to which existing theorems on such averages \cite{BreuerMajor1983,Hannan1970} can be readily applied. 
We further show that the asymptotics remain valid under modestly irregular sampling due to jitter or missing observations.
For the fixed-domain setting, the results extend to models that behave sufficiently like the power-law model at high frequencies such as the often used Mat\'ern model \cite{ZhuStein2006, WangLoh2011EJS, KaufmanShaby2017EJS}.
\end{abstract}



\end{frontmatter}

\section{Introduction}\label{sec:intro}

For describing the local behavior of a continuous, isotropic Gaussian random field, two key quantities are its roughness and scale.
The simplest model for the covariance structure that includes these quantities is the power-law Gaaussian random field, which arises in geophysics, climate reanalysis, and cosmic microwave background (CMB) temperature fields \cite{GneitingSchlather2004,Falconer2003,ChilesDelfiner2012,MarinucciPeccati2011}. 
The power-law random field is an intrinsic random function (IRF) with generalized covariance \cite{Matheron1973}
\begin{equation}
  K_{\phi_1,\phi_2}(h)=\phi_1\,\Gamma(-\phi_2)\,\|h\|^{2\phi_2},\qquad
  \phi_1>0,\ \phi_2>0
  \label{eq:powerlaw}
\end{equation}
when $\phi_2$ is not an integer.
For $\phi_2$ a positive integer,
\begin{equation}
  K_{\phi_1,\phi_2}(h)=\frac{2\phi_1(-1)^{\phi_2+1}}{\phi_2 !} \|h\|^{2\phi_2} \log\|h\|.
   \label{eq:power-integer}
 \end{equation}
In \(d=2\) dimensions, the associated continuous-space spectral density is
\begin{equation}
  f_0(\lambda)=c_{2,\phi_2}\,\phi_1\,\|\lambda\|^{-(2+2\phi_2)},\qquad
  c_{2,\phi_2}=\frac{2^{2\phi_2}\pi\,\Gamma(\phi_2)}{\Gamma(1+\phi_2)}.
  \label{eq:power-sd}
\end{equation}
The inclusion of the factor $\Gamma(-\phi_2)$ in the definition of $K_{\phi_1,\phi_2}$ in (\ref{eq:powerlaw}) may appear unnecessary, but if it is omitted from (\ref{eq:powerlaw}) then the corresponding formula for $c_{2,\phi_2}$ in (\ref{eq:power-sd}) equals 0 when $\phi_2$ is a positive integer.

\paragraph{Generalized covariances (GC–\(k\)).}
For an integer \(k\ge0\), a kernel \(K\) is a GC–\(k\) if every signed measure annihilating polynomials of total degree at most \( k\) produces a positive finite-variance linear functional of the field. The power-law form in \eqref{eq:powerlaw} and \eqref{eq:power-integer} is GC–\(k\) \emph{iff} \(\phi_2<k+1\).\footnote{For \(k=0\), the classical semivariogram equals \(-K\).}
When \(0<\phi_2<1\) the surface is fractal.
Larger $\phi_2$ yields smoother fields.
In particular, the Gaussian power-law field is $p$ times differentiable in any direction if and only if $\phi_2 > p$.

This work considers asymptotic properties of estimates of $\phi_1$ and $\phi_2$ when the random field is observed on a square grid.
Spatial statistical theory distinguishes two classical asymptotic regimes.
Under \emph{increasing-domain} (ID) asymptotics, the volume of the observation window
increases at the same rate as the number of observations, while under \emph{fixed-domain} (FD)
asymptotics the domain remains fixed and bounded as the number of observations increases.
FD asymptotics even for simple parametric models for the covariance function can be tricky because there can naturally be functions of the parameters that cannot be estimated consistently as the number of observations tends to infinity \cite{Ying1991,Zhang2004,Stein1999}.
However, for the power-law model, this problem does not occur as demonstrated by \cite{ZhuStein2002} when $\phi_2 < 1$.
In an unpublished thesis, \cite{ShinThesis} shows how both $\phi_1$ and $\phi_2$ can be consistently estimated for all positive $\phi_2$ using a spectral domain estimate, but does not quite give a result for the joint asymptotic distribution of the estimates.

In this work, for the power-law model, we first give a direct proof under increasing-domain (ID) asymptotics of consistency and joint asymptotic normality for a class of moment-based estimators of \(\phi_1\) and \(\phi_2\) using gridded observations. Exploiting the exact self-similar rescaling of the power-law model, we then transfer these limits to the fixed-domain (FD) setting without the need to reprove FD specific central limit theorems. This shortcut avoids several technicalities that arise in a direct FD analyses.

For context, many classical FD results for the Matérn class identify only a microergodic \cite{Stein1999} combination of parameters that assumes the roughness parameter is known \cite{Ying1991,Zhang2004}. In contrast, for the power law model, all parameters are consistently estimable in FD when \(0<\phi_2<1\) \cite{ZhuStein2002}. Our contribution provides joint CLTs for \((\hat\phi_1,\hat\phi_2)\) in both ID and FD regimes, with the FD limits obtained via the self-similar mapping. See \cite{Zhang2005} for a discussion of asymptotic regimes in spatial statistics. Related spectral/Whittle approaches primarily target tail index (roughness) estimation rather than joint scale–roughness inference; cf.\ \cite{Stein1995JASA,LimStein2008,ShinThesis}. Classical increasing-domain spectral results in lattices (e.g., \citep{Guyon1982, DahlhausKunsch1987}) assume that the (aliased) spectral density is uniformly bounded away from zero and infinity and sufficiently smooth on $(-\pi,\pi]^d)$, assumptions violated by filtered power–law fields with a spectral pole and limited differentiability at the origin. More recent edge/aliasing corrections such as the debiased spatial Whittle likelihood \citep{GuillauminEtAl2022} retain similarly regular spectral conditions and, therefore, still exclude intrinsic/power-law regimes; our analysis targets precisely this nonregular case.

\paragraph{Variation-based inference.}
Finite-difference \emph{quadratic variations} and wavelet methods
offer likelihood-free routes to roughness estimation.
In one dimension, quadratic variations of fractional Brownian motion
yield central limit theorems for the Hurst parameter
\cite{IstasLang1997,Coeurjolly2001}, and higher-order variations have
been analyzed in related contexts. Applications in finance and econometrics have exploited such methods
to quantify roughness in asset returns and volatility surfaces
\cite{GatheralJaissonRosenbaum2018,ElEuchRosenbaum2019}.
In higher dimensions, extensions include generalized quadratic variations for fractional Brownian sheets \cite{BiermeBonamiLeon2011-EJP}, anisotropic surface models, and wavelet-based estimators of local Hurst exponents under irregular designs. Most of these results develop distributional theory for the \emph{roughness} (Hurst/smoothness) parameter with the overall scale profiled out or treated as a nuisance, and their asymptotics are typically derived under long span/increasing-domain setups (fixed grid spacing with the window growing) rather than fixed-domain limits. Fixed-domain distributional results exist mainly for specific fractional models and still focus on roughness (e.g., \cite{ZhuStein2002}); joint scale–roughness FD theory in more than one dimension appears to be limited.

\paragraph{Our contribution.}
Similar to \cite{ZhuStein2002}, we propose a \emph{method of moments} for joint estimation of both scale and roughness parameters in two-dimensional power-law IRFs.
We choose a specific form for these estimates that simplifies the theory and yields reasonably efficient parameter estimates.
Our main results are:
\begin{enumerate}
\item \textbf{Joint distributional theory.}  
      We prove simultaneous consistency and central limit theorems for
      both parameters under ID and FD asymptotics, providing the first
      unified framework for joint inference in both regimes.
      This unified framework provides simple explanations for the following phenomena that occur for $N$ gridded observations under FD asymptotics: the rate of convergence for the estimated scale parameter is slightly slower than the standard $N^{-1/2}$ rate and the correlation between the estimates of scale and roughness tends to 1 as $N\to\infty$.
\item \textbf{Extension to smoother processes.}  
      By employing higher-order difference filters we extend the
      methodology beyond fractal surfaces to smoother
      processes with roughness exponents greater than one.
\item \textbf{Robustness.}  
      We establish robustness to practical deviations, including site
      deletion, Bernoulli thinning, spatial jitter, and in case of the fixed domain setting, under models that behave only locally like a power law.
\end{enumerate}
Together these contributions provide a theoretically
rigorous and computationally scalable alternative to likelihood-based methods, bridging the gap between classical roughness estimation and joint parameter inference for fractal and smoother Gaussian surfaces.

\section{Related Work}\label{sec:litreview}

\paragraph{One-dimensional models.}
Much of the early literature focused on one-dimensional fractional
Brownian motion (fBm), whose increments exhibit long-range dependence.
Quadratic variation and increment-based estimators of the Hurst
exponent were introduced by \cite{IstasLang1997} and extended in
\cite{Coeurjolly2001}, yielding central limit theorems under
increasing-domain asymptotics.  
Wavelet-based methods provide alternatives with multiscale localization
\cite{AbryVeitch1998}, and higher-order
variations have been analysed in related Gaussian process settings.  
Such approaches have been widely applied in finance and econometrics
to quantify rough volatility and long memory in asset returns
\cite{GatheralJaissonRosenbaum2018,ElEuchRosenbaum2019}.  
However, these contributions primarily address estimation of the
roughness parameter, treating scale as a nuisance and focusing on
time series data. In the fixed-domain (generally called ``high-frequency'' in the finance literature) setting, \cite{kawai2013a,kawai2013b} provide theory for the joint estimateion of roughness and scale for fBm-type models, complementing the above ID-focused literature.

\paragraph{Fractional Brownian surfaces and Gaussian fields.}
The two-dimensional analogue of fBm, fractional Brownian surface (fBs),
has been studied in both spatial statistics and probability.
By considering the FD asymptotic behavior of the Fisher information matrix under a periodic analog to the Mat\'ern model, \cite{Stein1999} appears to have been the first to notice that in any number of dimensions, the scale parameter estimate is likely to have slower than $N^{-1/2}$ convergence and correlation with the roughness parameter estimate tending to 1.
\cite{Chan2000} studied fixed-domain asymptotic properties of increment-based estimates of the roughness parameters.
\cite{ZhuStein2002} considered estimation of both roughness and scale for increment-based estimates similar to those considered in this work and gave FD asymptotic results for the marginal asymptotic distribution of both parameter estimates, but did not consider their joint distribution.  
Subsequent work established quadratic-variation central limit theorems
for Gaussian fields \cite{BiermeBonamiLeon2011-EJP}, Radon transform
methods for anisotropic fractional surfaces,
and wavelet-based estimation of local Hurst exponents on irregular
designs.
The results in this work, however, focus on the roughness parameter, with scale either fixed or regarded as a nuisance.

\paragraph{Likelihood-based and related approaches.}
For Matérn covariances, FD asymptotics imply only microergodic combinations are consistently estimable \cite{Ying1991,Ying1993,Zhang2004,Stein1999}. 
Tapered and compactly supported models can reduce computations for likelihood-based methods while maintaining asymptotically optimal estimation in some cases under FD asymptotics \cite{DuZhangMandrekar2009,WangLoh2011EJS,Gneiting2002,BevilacquaEtAl2019}, and the role of the range parameter for estimation/prediction is analysed in \cite{KaufmanShaby2013}. Additional perspectives include cross-validation under misspecification \cite{Bachoc2013}, the role of spatial sampling \cite{Bachoc2014JMA}, and FD asymptotics when the Gaussian process is observed with measurement error \cite{ChenSimpsonYing2000}. 

While powerful, likelihood-based approaches require heavy computation, whereas the moment estimators we study here are explicit and scale linearly in data size. Frequency-domain methods based on the (tapered/debiased) spatial Whittle likelihood reduce the cost of likelihood evaluation to $\mathcal O(N\log N)$ via FFTs on an $n\times n$ grid ($N=n^2$), and under a correctly specified parametric stationarity model, can attain near-efficient inference after addressing edge effects (e.g., by tapering and/or using the expected/debiased periodogram). See \cite{Whittle1954, DahlhausKunsch1987, GuillauminEtAl2022}.  By contrast, our two–scale quadratic-variation estimators are fully explicit and require only a single pass with fixed stencils, i.e., $\mathcal O(N)$ time and $\mathcal O(1)$ working memory, with no global transforms. They remain valid for intrinsically nonstationary IRF-$(m-1)$ fields after $m$'th-order differencing (Section~\ref{subsec:DiffM}) and for Matérn truth with the same local roughness. While Whittle-type procedures can be more statistically efficient under correctly specified stationary ID models, our approach trades a small amount of efficiency for linear complexity, locality, and robustness to intrinsic structure \cite{Whittle1954, DahlhausKunsch1987, GuillauminEtAl2022}.

For $\phi_2 < 1$, \cite{ZhuStein2002} describe a class of moment-based estimates of roughness and scale that include our estimates as a special case.
While this work gives FD asymptotic distributions for each parameter estimate for fractional Brownian surfaces, it does not consider the joint distribution of the two estimates.

Using a smoothed periodogram, \cite{ShinThesis} develops estimators for the high-frequency behavior of the spectral density under parametric families that essentially include a scale parameter and a \emph{spectral-tail} (high-frequency) exponent that corresponds to $\phi_2$ here.
Unlike \cite{ZhuStein2002}, this work treats all $\phi_2 > 0$.
Assuming stationarity, \cite{ShinThesis} obtains FD asymptotic results for both the roughness and scale parameters, but does not give explicit results on their joint asymptotic distribution such as those provided in Section \ref{subsec:FD-CLTs-1}.

By contrast, our self–affine, differencing based analysis yields FD \emph{joint} asymptotics for both roughness and scale estimates.
Furthermore, by exploiting the connection between the FD and the ID settings under the power-law model, the source of the non-standard asymptotics in the FD setting becomes transparent.
Our work also provides theoretical results showing when common sampling perturbations leave asymptotic distributions unchanged and shows that the same FD asymptotic results hold if the true covariance function is Mat\'ern rather than power-law.

\paragraph{Quadratic forms and CLTs.}
Finally, our work connects to the general theory of quadratic forms of
Gaussian fields.  
Classical results such as \cite{BreuerMajor1983} established central
limit theorems for nonlinear functionals of Gaussian sequences, while
\cite{AvramLeonenkoSakhno2010-ESAIM} proved
Szegő-type limit theorems and central limit theorems for quadratic
forms of stationary fields.  
These foundational results underpin the asymptotics of
variation-based estimators.  
We leverage them to establish joint central limit theorems for both
roughness and scale parameters, filling a gap left by prior work that
estimated only one at a time.

\section{Methodology: Bilinear Differences and Moment Estimators}
\label{sec:method}

\paragraph{Sampling and aliasing.}
Sampling the random field on the integer lattice \(\mathbb Z^2\) aliases the spectrum to the torus \(\mathbb T^2=(-\pi,\pi]^2\): 
\[
  f_X(\lambda)=f_0^{\mathrm{lat}}(\lambda)=\sum_{m\in\mathbb Z^2} f_0(\lambda+2\pi m)
\]
for $f_0$ as defined in (\ref{eq:power-sd}).
For convenience, define 
  \( \mu(\dd\lambda)=(2\pi)^{-2}\dd\lambda \).

\paragraph{Two bilinear differences (first order).}
Assume we observe $X$ on the set \(\Lambda_n=\{0,\dots,n-1\}^2\).
Define \(e_1=(1,0)\), \(e_2=(0,1)\). We use first-order bilinear product differences at steps \(r\in\{1,2\}\):
\begin{equation}
  (D^{(1)}_{[r]}X)_t
  = X_t - X_{t+re_1} - X_{t+re_2} + X_{t+re_1+re_2},
  \qquad r\in\{1,2\},
  \label{eq:bilinear-def}
\end{equation}
with symbols \begin{equation}
  g^{[r]}_1(\lambda)=(1-e^{ir\lambda_1})(1-e^{ir\lambda_2}),
  \qquad |g^{[r]}_1(\lambda)|^2\le 16 .
  \label{eq:bilinear-symbol}
\end{equation}
We have
\begin{lemma}[Expectation for $D^{(1)}_{[r]}$]\label{lem:a1}
If $0<\phi_2<1$, then
\[
   \E\bigl[(D^{(1)}_{[r]}X_t)^2\bigr]
   = \phi_1\,2^{2(r-1)\phi_2} a_1(\phi_2) \;\mbox{ with }\;
   a_1(\phi_2)=|\Gamma(-\phi_2)|(8-4\cdot 2^{\phi_2}).
\]
\end{lemma}

\begin{proof}
See Appendix \ref{proofs}.
\end{proof}
See \S\ref{sec:HigherOrder} for the higher-order analogue when \(\phi_2\ge 1\)).

\paragraph{Quadratic variations and two-scale identity.}
Our estimators of $\phi_1$ and $\phi_2$ are based on averages of observed squares of these bilinear differences:
\begin{equation}
  Q^{(1)}_r
  =\frac{1}{|\Lambda_{n-2r}|}\sum_{t\in\Lambda_{n-2r}}
    \bigl(D^{(1)}_{[r]}X_t\bigr)^2,
  \label{eq:Q-def}
\end{equation}
where, $|A|$ is the size of a finite set $A$.
Defining $q^{(1)}_r=\E[Q^{(1)}_r]$, we have
\begin{equation}
  q^{(1)}_r
  = \int_{\mathbb T^2} |g^{[r]}_1(\lambda)|^2 f_X(\lambda)\,\mu(\dd\lambda)
  \; \mbox{ and } \;
  \frac{q^{(1)}_2}{q^{(1)}_1}=2^{\,2\phi_2}.
  \label{eq:Q-mean-ratio}
\end{equation}
This result suggests the following moment-based estimators:
\begin{equation}
  \widehat\phi_2=\tfrac12\log_2\!\left(\frac{Q^{(1)}_2}{Q^{(1)}_1}\right),
  \qquad
  \widehat\phi_1=\frac{Q^{(1)}_1}{a_1(\widehat\phi_2)} .
  \label{eq:phi-hats-m1}
\end{equation}

\paragraph{Notation.}
For two nonnegative functions $a(\cdot)$ and $b(\cdot)$ we write 
$a \lesssim b$ if there exists a constant $C>0$, independent of the 
argument, such that $a \le Cb$.  
We write $a \gtrsim b$ if $a \ge c b$ for some $c>0$, and 
$a \asymp b$ if both $a \lesssim b$ and $a \gtrsim b$ hold.

\paragraph{Integrability and the $D^{(1)}_{[2]}$–from–$D^{(1)}_{[1]}$ recovery.}
To apply the quadratic–form CLT of \cite[Thm.~2.2]{AvramLeonenkoSakhno2010-ESAIM}
we need (i) $f_{D^{(1)}_{[1]}X}\in L^2(\T^2)$ and (ii) bounded quadratic–form
symbols. Crucially, \(D^{(1)}_{[2]}X\) is a local linear transform of \(D^{(1)}_{[1]}X\).

\medskip\noindent
\textit{Spatial recovery of \(D^{(1)}_{[2]}X\) from \(D^{(1)}_{[1]}X\).}
\[
 g^{[2]}_1(\lambda)=g^{[1]}_1(\lambda)\,h(\lambda),\quad
 h(\lambda)=(1+e^{i\lambda_1})(1+e^{i\lambda_2}),
\]
so \(D^{(1)}_{[2]} = H\!\circ D^{(1)}_{[1]}\), where \(H\) is the four–point block–sum operator
\[
   (HZ)_t \;=\; Z_t+Z_{t+e_1}+Z_{t+e_2}+Z_{t+e_1+e_2}.
\]
Hence
\[
  (D^{(1)}_{[2]}X)_t \;=\; \bigl(H\,D^{(1)}_{[1]}X\bigr)_t,
\]
i.e. \(D^{(1)}_{[2]}X_t\) is a bounded linear combination of nearby
\(D^{(1)}_{[1]}X\) values. 
It is this fact that allows us to obtain asymptotic bivariate normality for $Q_1^{(1)}$ and $Q_2^{(1)}$ by applying \cite{AvramLeonenkoSakhno2010-ESAIM}.
We expect that similar results apply to a broad range of filters \cite{ZhuStein2002}, but would not follow from \cite{AvramLeonenkoSakhno2010-ESAIM}.
In the frequency domain,
\[
  B(\lambda):=\frac{|g^{[2]}_1(\lambda)|^2}{|g^{[1]}_1(\lambda)|^2}
  =|h(\lambda)|^2
  =16\cos^2\!\left(\frac{\lambda_1}{2}\right)\cos^2\!\left(\frac{\lambda_2}{2}\right)
  \le 16.
\]
Thus both $Q^{(1)}_1$ and $Q^{(1)}_2$ are quadratic forms of the \emph{same}
stationary field \(D^{(1)}_{[1]}X\) with bounded symbols
\[
  b_1(\lambda)\equiv 1,\qquad b_2(\lambda)=B(\lambda)\in L^\infty(\T^2).
\]

\medskip\noindent
\textit{Behavior near the origin.}
As $\theta\to0$, $|1-e^{i\theta}|^2=2(1-\cos\theta)\sim \theta^2$, hence
for small $\|\lambda\|$,
$|g^{[r]}_1(\lambda)|^2\lesssim r^4\|\lambda\|^4$, and away from the axes
$|g^{[r]}_1(\lambda)|^2\asymp r^4\|\lambda\|^4$.
For an IRF–0 power–law field in $d=2$,
$f_0(\lambda)\asymp \|\lambda\|^{-(2+2\phi_2)}$ as $\lambda\to0$, and the
aliasing sum is dominated by $k=0$ near $0$, so
$f_X(\lambda)=f_0^{\lat}(\lambda)\asymp f_0(\lambda)$ there. Consequently,
\[
  f_{D^{(1)}_{[r]}X}(\lambda) \;\lesssim\; \|\lambda\|^{\,2-2\phi_2}\qquad (\lambda\to0).
\]

\medskip\noindent
\textit{Square integrability.}
Writing $s=\|\lambda\|$, we have $f_{D^{(1)}_{[r]}X}(\lambda)^2 \lesssim s^{\,4-4\phi_2}$ and
\[
  \int_{\|\lambda\|\le \varepsilon} f_{D^{(1)}_{[r]}X}(\lambda)^2 \,\dd\lambda
   \;\lesssim\; \int_0^\varepsilon s^{\,4-4\phi_2}\, s\,\dd s
   = \int_0^\varepsilon s^{\,5-4\phi_2}\,\dd s,
\]
which converges at $0$ if $\phi_2<\tfrac32$.  
Hence $f_{D^{(1)}_{[r]}X}\in L^2(\T^2)$ for $r=1,2$ whenever $0<\phi_2<\tfrac32$.
\medskip

\textit{Boundedness of filter symbols.}
Since $|1-e^{ir\theta}|\le 2$ for all $\theta$, both $\big|g^{[1]}_1\big|^2$ and
$\big|g^{[2]}_1\big|^2$ are bounded on $\T^2$, and $b_2=B\le16$.  
Thus $g^{[r]}_1\in L^\infty(\T^2)$ and
$b_{\alpha,\beta}=\alpha b_1+\beta b_2\in L^\infty(\T^2)$ for any
scalars $(\alpha,\beta)$.

\begin{lemma}[Integrability of the filtered spectrum]\label{lem:Integrability}
Let $X$ be a two–dimensional IRF–0 power–law field with $0<\phi_2<1$.
Then $f_{D^{(1)}_{[1]}X}(\lambda)=|g^{[1]}_1(\lambda)|^2 f_X(\lambda)\in L^2(\T^2)$
and $g^{[1]}_1\in L^\infty(\T^2)$.  
Therefore the assumptions of \cite[Thm.~2.2]{AvramLeonenkoSakhno2010-ESAIM} hold,
and any linear combination $\alpha Q^{(1)}_1+\beta Q^{(1)}_2$ can be treated as a quadratic form of
the stationary field $D^{(1)}_{[1]}X$.
\end{lemma}

\section{Central limit theorems under increasing–domain asymptotics}
\label{sec:IDCLT}

We study the asymptotic distribution of the quadratic–variation
statistics as the observation window expands. Throughout this section
assume $0<\phi_2<1$, so the filtered spectrum
$f_{D^{(1)}_{[r]}X}\in L^2(\T^2)$ and the filter symbols are bounded
(Lemma~\ref{lem:Integrability}). Hence the quadratic–form central limit
theorem of \cite[Theorem~2.2]{AvramLeonenkoSakhno2010-ESAIM} applies with
$(p_f,p_b)=(2,\infty)$.

\subsection{Univariate CLT for linear combinations}
\label{sec:UniCLT}

Consider linear combinations of the quadratic–variation statistics
\[
   L_{\alpha,\beta}=\alpha Q^{(1)}_{1}+\beta Q^{(1)}_{2},
   \qquad (\alpha,\beta)\in\R^2\setminus\{(0,0)\},
\]
so that $\E[L_{\alpha,\beta} = q_{\alpha,\beta}=\alpha q^{(1)}_{1}+\beta q^{(1)}_{2}$.

It is convenient to view both $Q^{(1)}_{1}$ and $Q^{(1)}_{2}$ as quadratic forms of the same stationary
increment field $D^{(1)}_{[1]}X$, whose (lattice) spectrum is
$f_{D^{(1)}_{[1]}X}(\lambda)=|g^{[1]}_1(\lambda)|^2 f_X(\lambda)$. Then
\[
   q^{(1)}_{1}=\int_{\T^2} f_{D^{(1)}_{[1]}X}(\lambda)\,\mu(\dd\lambda),
   \qquad
   q^{(1)}_{2}=\int_{\T^2} \frac{|g^{[2]}_1(\lambda)|^2}{|g^{[1]}_1(\lambda)|^2}\,
              f_{D^{(1)}_{[1]}X}(\lambda)\,\mu(\dd\lambda).
\]
Define the bounded symbols
\[
   b_1(\lambda)\equiv 1,
   \qquad
   b_2(\lambda)=\frac{|g^{[2]}_1(\lambda)|^2}{|g^{[1]}_1(\lambda)|^2}
   =16\cos^2\!\Bigl(\frac{\lambda_1}{2}\Bigr)\cos^2\!\Bigl(\frac{\lambda_2}{2}\Bigr)\le 16,
\]
and set $b_{\alpha,\beta}=\alpha b_1+\beta b_2$.

\begin{lemma}[Boundary remainder is negligible at $\sqrt{N}$ scale]\label{lem:trim-slutsky}
Let
\[
  Q_{[1]}^\circ=\frac{1}{|\Lambda_{n-4}|}\sum_{t\in\Lambda_{n-4}}\bigl(D^{(1)}_{[1]}X_t\bigr)^2,
\]
and write $N=n^2$. If $0<\phi_2<1$, then 
\[
  \E Q^{(1)}_{1}=\E Q_{[1]}^\circ
  \qquad\text{and}\qquad
  \Var\Bigl(\sqrt{N}\,\bigl[Q^{(1)}_{1}-Q_{[1]}^\circ\bigr]\Bigr)\longrightarrow 0,
\]
hence $\sqrt{N}\,[Q^{(1)}_{1}-Q_{[1]}^\circ]\xrightarrow{L^2}0$ (in particular $\xrightarrow{p}0$).
\end{lemma}

\begin{proof}
    See Appendix \ref{proofs}
\end{proof}

\noindent
Define the trimmed linear combination
\[
  L_{\alpha,\beta}^{\circ}=\alpha Q_{[1]}^\circ+\beta Q^{(1)}_{2},
\]
so both terms are quadratic forms of $D^{(1)}_{[1]}X$ over the \emph{same} index set
$\Lambda_{n-4}$ with symbol $b_{\alpha,\beta}$. 

\begin{proposition}[Quadratic–form CLT under increasing domain]
\label{thm:UniCLT}
Assume $0<\phi_2<1$ so that $f_{D^{(1)}_{[1]}X}\in L^2(\T^2)$ and $b_{\alpha,\beta}\in L^\infty(\T^2)$.
For any fixed $(\alpha,\beta)\neq(0,0)$,
\[
   \sqrt{N}\,\bigl(L_{\alpha,\beta}-q_{\alpha,\beta}\bigr)
   \;\Rightarrow\;
   \mathcal N\!\bigl(0,\sigma_{\alpha,\beta}^2\bigr),
\]
where
\[
   \sigma_{\alpha,\beta}^2
   \;=\;
   2\int_{\T^2} b_{\alpha,\beta}(\lambda)^2\, f_{D^{(1)}_{[1]}X}(\lambda)^2\,\mu(\dd\lambda).
\]
\end{proposition}

\begin{proof}
By Lemma~\ref{lem:Integrability}, $f_{D^{(1)}_{[1]}X}\in L^2(\T^2)$ and $b_{\alpha,\beta}\in L^\infty(\T^2)$,
so \cite[Thm.~2.2]{AvramLeonenkoSakhno2010-ESAIM} yields
\[
  \sqrt{N}\,\bigl(L_{\alpha,\beta}^{\circ}-\E L_{\alpha,\beta}^{\circ}\bigr)
  \Rightarrow \mathcal N(0,\sigma_{\alpha,\beta}^2).
\]
By Lemma~\ref{lem:trim-slutsky},
$\sqrt{N}\,\bigl[L_{\alpha,\beta}-L_{\alpha,\beta}^{\circ}\bigr]\xrightarrow{p}0$.
Slutsky’s theorem gives the stated limit for $L_{\alpha,\beta}$.
\end{proof}

\begin{remark*}[Normalizing by $N$ vs.\ $|\Lambda_{n-2r}|$]
Since $|\Lambda_{n-2r}|/N=(1-2r/n)^2\to1$ as $n\to\infty$, replacing $|\Lambda_{n-2r}|$ by $N$
only rescales by a factor tending to one. In particular, multiplying by $\sqrt{N}$ (or
$\sqrt{|\Lambda_{n-2r}|}$) yields the same $\sqrt{N}$ limits and asymptotic covariance.
\end{remark*}
The preceding proposition yields joint asymptotic normality by the
Cramér–Wold device.

\begin{theorem}[Joint CLT for $(Q^{(1)}_{1},Q^{(1)}_{2})$]\label{thm:Joint}
Let
\[
  b_1(\lambda)\equiv 1,\qquad
  b_2(\lambda)=\frac{|g^{[2]}_1(\lambda)|^2}{|g^{[1]}_1(\lambda)|^2}
  =16\cos^2\!\left(\frac{\lambda_1}{2}\right)\cos^2\!\left(\frac{\lambda_2}{2}\right)\le 16,
\]
and $f_{D^{(1)}_{[1]}X}(\lambda)=|g^{[1]}_1(\lambda)|^2 f_X(\lambda)$. Define
\begin{equation}
  \Sigma_{\ell m}
  =2\int_{\mathbb T^2} b_\ell(\lambda)b_m(\lambda)
    \,f_{D^{(1)}_{[1]}X}(\lambda)^2\,\mu(\dd\lambda),
  \quad b_1\equiv1,\quad
  b_2(\lambda)=16\cos^2\!\tfrac{\lambda_1}{2}\cos^2\!\tfrac{\lambda_2}{2}.
  \label{eq:Sigma-m1}
\end{equation}

Assume $0<\phi_2<1$. Then
\[
  \sqrt{N}\!
  \begin{pmatrix}
     Q^{(1)}_{1}-q^{(1)}_{1}\\[2pt]
     Q^{(1)}_{2}-q^{(1)}_{2}
  \end{pmatrix}
  \Rightarrow \mathcal N(\bm 0,\Sigma).
\]
\end{theorem}

\subsection{Method-of-moments estimators and delta–method CLT}

\begin{theorem}[Consistency and joint asymptotic normality]
\label{thm:Delta}
For $0<\phi_2<1$, we have:
\begin{itemize}
\item \textbf{Consistency:}\;
$(\widehat\phi_1,\widehat\phi_2)\xrightarrow{p}(\phi_1,\phi_2)$.

\item \textbf{Asymptotic normality:}\;
Define
\begin{equation}
  h(y_1,y_2)
  =\Biggl(
     \frac{y_1}{a_1\!\left(\tfrac12\log_2(y_2/y_1)\right)},
     \ \tfrac12\log_2(y_2/y_1)
    \Biggr).
  \label{eq:delta-map}
\end{equation}

and let $J=\nabla h(q^{(1)}_{1},q^{(1)}_{2})$ be the Jacobian at $(q^{(1)}_{1},q^{(1)}_{2})$. Then
\[
  \sqrt{N}
  \begin{pmatrix}
     \widehat\phi_1-\phi_1\\[3pt]
     \widehat\phi_2-\phi_2
  \end{pmatrix}
  \;\Rightarrow\;
  \mathcal N\!\bigl(\bm 0,\,J\,\Sigma\,J^\top\bigr),
\]
where $\Sigma$ is as in Theorem~\ref{thm:Joint}.
\end{itemize}
\end{theorem}

\begin{proof}
By Theorem~\ref{thm:Joint}, $\sqrt{N}\!\left((Q^{(1)}_{1},Q^{(1)}_{2})-(q^{(1)}_{1},q^{(1)}_{2})\right)\Rightarrow\mathcal N(\bm 0,\Sigma)$.
The map $h$ is $C^1$ on $(0,\infty)^2$ because $a_1(\phi_2)>0$ and smooth on $(0,1)$, and
$\log_2$ is smooth on $(0,\infty)$. Hence the multivariate delta method
gives the stated limit. Consistency follows by the continuous mapping
theorem applied to $Q^{(1)}_{j}\to q^{(1)}_{j}$ in probability and the continuity
of $h$.
\end{proof}

\begin{remark*}[Explicit Jacobian entries]
Let $A(\phi)=a_1(\phi)$ and put $s=(2\ln 2)^{-1}$.  
With $r=y_2/y_1$,
\[
  h_2(y_1,y_2)=\tfrac12\log_2 r
  \;\Rightarrow\;
  \frac{\partial h_2}{\partial y_1}=-\frac{s}{y_1},\qquad
  \frac{\partial h_2}{\partial y_2}=\frac{s}{y_2}.
\]
For $h_1(y_1,y_2)=y_1/A(h_2(y_1,y_2))$,
\[
  \frac{\partial h_1}{\partial y_1}
    = \frac{1}{A(h_2)} + \frac{A'(h_2)}{A(h_2)^2}\,s,
  \qquad
  \frac{\partial h_1}{\partial y_2}
    = -\,\frac{A'(y_2)}{A(y_2)^2}\,s\,\frac{y_1}{y_2}.
\]
Evaluating at $(y_1,y_2)=(q^{(1)}_{1},q^{(1)}_{2})$ (so $h_2=\phi_2$, $A=a_1(\phi_2)$, $A'=a_1'(\phi_2)$) gives
\[
  J=
  \begin{pmatrix}
    A^{-1} + s\,A'/A^{2} & -\,s\,A'/A^{2}\,(q^{(1)}_{1}/q^{(1)}_{2})\\[4pt]
    -\,s/q^{(1)}_{1} & s/q^{(1)}_{2}
  \end{pmatrix},
\]
and $q^{(1)}_{1}/q^{(1)}_{2}=2^{-2\phi_2}$.
\end{remark*}

\begin{corollary}[Ratio form of the joint CLT]
Under the conditions of Theorem~\ref{thm:Delta},
\[
  \sqrt{N}
  \begin{pmatrix}
    \widehat\phi_1/\phi_1 - 1\\[2pt]
    \widehat\phi_2 - \phi_2
  \end{pmatrix}
  \ \Rightarrow\
  \mathcal N\!\bigl(\bm 0,\ K\,\Sigma_h\,K^\top\bigr),
  \qquad
  K=\begin{pmatrix}1/\phi_1&0\\[2pt]0&1\end{pmatrix},
  \ \ \Sigma_h:=J\,\Sigma\,J^\top.
\]
\end{corollary}

\section{From increasing to fixed domain}\label{sec:IDtoFD}

For $0<\phi_2<1$,  we transfer the ID ratio–CLTs of §\ref{sec:IDCLT} to the fixed–domain (FD).

\subsection{Re–indexing and amplitude rescaling}\label{subsec:FD-reindex}
Let $N=|\Lambda_n|=n^2$ and define the FD re–indexed field $X'(t):=X(t/n)$ for $t\in\Lambda_n\subset\mathbb Z^2$.
Using the power–law GC–0 form, the increment variance satisfies
\begin{equation}\label{eq:fd-increments}
  \E\!\left[(X'(t+h)-X'(t))^2\right]
  \;=\;2\,\gamma_{\phi_1,\phi_2}(h/n)
  \;=\;2\,\gamma_{\tau_1,\phi_2}(h),\qquad h\in\mathbb Z^2,
\end{equation}
with the FD parameters
\begin{equation}\label{eq:fd-params}
\tau_1:=\phi_1\,N^{-\phi_2},\qquad \tau_2:=\phi_2\
\end{equation}
Thus $\{X(t/n):t\in\Lambda_n\}\stackrel{d}{=}\{X'(t):t\in\Lambda_n\}$, i.e. a unit–lattice sample with parameters $(\tau_1,\tau_2)$.
Equivalently, the amplitude–rescaled field $\widetilde X(t):=n^{\phi_2}X'(t)$ has the same law on $\Lambda_n$ as the ID field with $(\phi_1,\phi_2)$.

\subsection{Exact FD–ID mapping for the first–order estimators}\label{subsec:FD-map-1}
Let $D_{[r]}^{(1)}$, $Q_r^{(1)}$ be as in previous sections, and let
$(\widehat\phi_1,\widehat\phi_2)$ (ID) and $(\widehat\tau_1,\widehat\tau_2)$ (FD) be the MoM estimators defined in \eqref{eq:phi-hats-m1}.
Since $D_{[r]}^{(1)}$ is linear and $Q_r^{(1)}$ are homogeneous of degree two, amplitude rescaling multiplies by the same deterministic factor; therefore
\begin{equation}\label{eq:fd-id-ratio-1}
  \left(\frac{\widehat\tau_1}{\tau_1},\ \widehat\tau_2\right)\ \stackrel{d}{=}\ \left(\frac{\widehat\phi_1}{\phi_1},\ \widehat\phi_2\right)
\end{equation}

\subsection{FD ratio–CLTs}\label{subsec:FD-CLTs-1}
Let $\Sigma$ be the covariance in Theorem~\ref{thm:Joint} and $K$ the Jacobian from Theorem~\ref{thm:Delta}; set $\Omega:=K\Sigma K^\top$.
From \eqref{eq:fd-id-ratio-1} and the ID ratio–CLT,
\begin{equation}\label{eq:FD1-CLT}
   \sqrt{N}
  \begin{pmatrix}
    \widehat\tau_1/\tau_1-1\\[2pt]
    \widehat\tau_2-\tau_2
  \end{pmatrix}
  \ \Rightarrow\ \mathcal N\!\big(\bm 0,\ \Omega\big).  
\end{equation}

\paragraph{Equivalent plug–in/log forms.}
Since $\tau_1=\phi_1N^{-\phi_2}$, write $U_N:=\widehat\tau_1/\tau_1-1$ and $V_N:=\widehat\tau_2-\tau_2$, and define
\[
  W_N:=\frac{N^{\widehat\tau_2}\widehat\tau_1}{\phi_1}-1=(1+U_N)\exp\!\big((\log N)V_N\big)-1.
\]
A second–order expansion together with \eqref{eq:FD1-CLT} yields
\begin{equation}\label{eq:FD1-plugin}
  \sqrt{N}
  \begin{pmatrix}
    W_N-(\log N)(\widehat\tau_2-\tau_2)\\[4pt]
    \widehat\tau_2-\tau_2
  \end{pmatrix}
  \Rightarrow \mathcal N\!\big(\bm 0,\ \Omega\big),\qquad
  \mathrm{Corr}\!\big(\log(N^{\widehat\tau_2}\widehat\tau_1),\ \widehat\tau_2\big)\to 1.
\end{equation}
Equivalently, with
\[
  A_N:=\begin{pmatrix}1&-\log N\\[2pt]0&1\end{pmatrix},
\]
we obtain the stabilized log form
\begin{equation}\label{eq:FD1-matrix}
  \sqrt{N}\,A_N
  \begin{pmatrix}
    \log(\widehat\tau_1/\tau_1)\\[4pt]
    \widehat\tau_2-\tau_2
  \end{pmatrix}
  \Rightarrow \mathcal N\!\big(\bm 0,\ \Omega\big).
\end{equation}

\section{Higher–order differences and extension to smoother processes}
\label{sec:HigherOrder}

We extend the methodology and asymptotic results to arbitrary
$\phi_{2}>0$ by using higher–order bilinear differences. When
$0<\phi_{2}<1$, the first–order bilinear difference suffices.
For general $\phi_{2}\in(k,k+1)$ with $k\in\mathbb N$, we take an
$m$th–order bilinear product difference with $m\ge k+1$.

\subsection{m-th–order bilinear differences and their spectra}
\label{subsec:DiffM}

For an integer $m\ge1$, define the $m$th–order 2D bilinear difference
\begin{equation}
  (D^{(m)}X)_{t}
  = \sum_{a_1=0}^{m}\sum_{a_2=0}^{m}
      c^{(m)}_{(a_1,a_2)}\,
      X_{t+a_{1}e_{1}+a_{2}e_{2}},
  \qquad
  c^{(m)}_{(a_1,a_2)} := (-1)^{a_{1}+a_{2}}
      \binom{m}{a_{1}}\binom{m}{a_{2}} .
  \label{eq:dm-def}
\end{equation}
and its (angular–frequency) symbol and modulus
\begin{equation}
  g_{m}(\lambda)=(1-e^{i\lambda_{1}})^{m}(1-e^{i\lambda_{2}})^{m},
  \qquad
  |g_{m}(\lambda)|^{2}=4^{2m}
    \prod_{i=1}^{2}\sin^{2m}\Bigl(\frac{\lambda_{i}}{2}\Bigr).
  \label{eq:gm-symbol}
\end{equation}

For step $j\in\{1,2\}$ we denote by $D^{(m)}_{[j]}$ the same stencil
applied on squares of side $j$, with symbol
$g_{m}^{[j]}(\lambda)=(1-e^{ij\lambda_{1}})^{m}(1-e^{ij\lambda_{2}})^{m}$.

\begin{lemma}[Stationarity after differencing]\label{lem:stat_after_diff}
Let $k=\lfloor \phi_2\rfloor$ and $m\ge k+1$. For $j\in\{1,2\}$ the field
$D^{(m)}_{[j]}X$ is centered, second–order stationary on $\Z^2$ and
\[
  \Cov\!\bigl((D^{(m)}_{[j]}X)(t),\,(D^{(m)}_{[j]}X)(t+u)\bigr)
  \,=\,
  \sum_{a,b\in\{0,\dots,m\}^2} c_a^{(m)}c_b^{(m)}\,
    K_{\phi_1,\phi_2}\!\bigl(j(u+b-a)\bigr).
\]
\end{lemma}

\begin{proof}
Refer to \ref{proofs}
\end{proof}

\begin{lemma}[Spectra in the continuum and on the torus]\label{lem:spectrum}
Let $m\ge \lfloor \phi_{2}\rfloor+1$. Then:
\begin{enumerate}
\item (Continuum) The spectral density of $D^{(m)}_{[1]}X$ on $\mathbb R^2$ is
\[
   f^{\mathrm{cont}}_{D^{(m)}_{[1]}X}(\omega)
   \;=\;|g_{m}(\omega)|^{2}\,f_{0}^{\mathrm{cont}}(\omega)
   \;=\; |g_m(\omega)|^2\,C_{\phi_{2}}\phi_{1}\,\|\omega\|^{-(2+2\phi_{2})}.
\]
\item (Torus / aliased) On $\T^{2}$ the aliased spectrum is
\[
  f_{D^{(m)}_{[1]}X}(\lambda)
  :=\sum_{k\in\Z^{2}} f^{\mathrm{cont}}_{D^{(m)}_{[1]}X}(\lambda+2\pi k)
  \;=\; |g_{m}(\lambda)|^{2}\,f_{0}^{\lat}(\lambda).
\]
\end{enumerate}
\end{lemma}

\begin{proof}
Refer to \ref{proofs}
\end{proof}

\begin{lemma}[Square–integrability on $\T^2$]\label{lem:L2m}
Let $\phi_{2}\in(k,k+1)$ and $m\ge k+1$. Then
\(f_{D^{(m)}_{[1]}X}\in L^{2}(\T^{2})\).
\end{lemma}
\begin{proof}
Refer to \ref{proofs}
\end{proof}

\paragraph{Two–scale factorization.}
For each $m\ge1$,
\[
  g_m^{[2]}(\lambda)= g_m^{[1]}(\lambda)\,h_m(\lambda),
  \qquad
  h_m(\lambda) := (1+e^{i\lambda_1})^{m}(1+e^{i\lambda_2})^{m}.
\]
Let $H^{(m)}$ be the \((m{+}1)^2\)-point block–sum operator
\[
  (H^{(m)}Z)_t
  := \sum_{\alpha\in\{0,\dots,m\}^2}
      \binom{m}{\alpha_1}\binom{m}{\alpha_2}\,Z_{t+\alpha_1 e_1+\alpha_2 e_2}.
\]
Then $D^{(m)}_{[2]} = H^{(m)}\!\circ D^{(m)}_{[1]}$, i.e.
$D^{(m)}_{[2]}X_t=(H^{(m)}\,D^{(m)}_{[1]}X)_t$. Moreover,
\[
  B_m(\lambda):=\frac{|g_m^{[2]}(\lambda)|^2}{|g_m^{[1]}(\lambda)|^2}
  = |h_m(\lambda)|^2
  = \prod_{i=1}^2 \bigl(2\cos(\lambda_i/2)\bigr)^{2m}
  \le 4^{2m}\in L^\infty(\T^2).
\]
Thus $Q^{(m)}_1$ and $Q^{(m)}_2$ are quadratic forms of the \emph{same} stationary field
$D^{(m)}_{[1]}X$ with bounded symbols $b_1\equiv1$ and $b_2=B_m$.

\subsection{Two–scale quadratic variations and the exact scale ratio}
\label{subsec:TwoScaleM}

Let \(F(\lambda)=f_{D^{(m)}_{[1]}X}(\lambda)=|g_m(\lambda)|^2 f_0^{\lat}(\lambda)\).
With $B_m$ as above, $f_{D^{(m)}_{[2]}X}=B_m\,F$. For $j\in\{1,2\}$,
\[
   Q^{(m)}_{j}
   = \frac{1}{|\Lambda_{n-jm}|}\sum_{t\in\Lambda_{n-jm}}
     \bigl(D^{(m)}_{[j]}X_t\bigr)^2,
\]
so \(Q^{(m)}_{1}\) and \(Q^{(m)}_{2}\) are quadratic forms of the \emph{same}
stationary field \(D^{(m)}_{[1]}X\) with bounded kernels. Denote expectations
by \(q^{(m)}_{j}=\E[Q^{(m)}_{j}]\).

\paragraph{Explicit mean and the function \(a_m(\phi_2)\).}
By the allowable–measure variance identity for filtered IRF fields,
\[
  q^{(m)}_{j}
  = -\,\phi_1\,|\Gamma(-\phi_2)|\sum_{a,b} c^{(m)}_a c^{(m)}_b\,
       \|\,j(b-a)\,\|^{2\phi_2}
  \;=\; \phi_1\,a_m(\phi_2)\,j^{2\phi_2},
\]
where the scale function \(a_m:(0,\infty)\to(0,\infty)\) is given explicitly by
\[
  a_m(\phi_2)
  := -\,|\Gamma(-\phi_2)|\sum_{a,b} c^{(m)}_a c^{(m)}_b\,\|b-a\|^{2\phi_2}
  \;=\; -\,2\,|\Gamma(-\phi_2)|\!\!\sum_{a<b} c^{(m)}_a c^{(m)}_b\,\|b-a\|^{2\phi_2}.
\]
In particular, for the bilinear stencil ($m{=}1$) one recovers
\[
  a_1(\phi_2)=|\Gamma(-\phi_2)|\bigl(8-4\,2^{\phi_2}\bigr),
\]
since four edge pairs contribute \(-1\) at distance \(1\) and two diagonal
pairs contribute \(+1\) at distance \(\sqrt2\).

\begin{proposition}[Exact expectation ratio]\label{prop:RatioM}
For $\phi_{2}>0$ and $m\ge \lfloor\phi_{2}\rfloor+1$,
\(
  q^{(m)}_{2}/q^{(m)}_{1}=2^{\,2\phi_{2}}.
\)
Consequently,
\[
   \widehat\phi_{2}
      =\tfrac12\log_{2}\!\left(\frac{Q^{(m)}_{2}}{Q^{(m)}_{1}}\right),
   \qquad
   \widehat\phi_{1}
      =\frac{Q^{(m)}_{1}}{a_{m}(\widehat\phi_{2})}.
\]
\end{proposition}
\begin{proof}
Refer to \ref{proofs}
\end{proof}

\paragraph{Common index set and negligible boundary.}
Let $N=n^2$ and write
\[
  Q^{(m),\circ}_j
   := \frac{1}{|\Lambda_{n-2m}|}\sum_{t\in\Lambda_{n-2m}} \bigl(D^{(m)}_{[j]}X_t\bigr)^2,
   \qquad j=1,2.
\]
As in Lemma~\ref{lem:trim-slutsky} for $m{=}1$, the boundary layer has $O(n)$ sites
(with width $O(m)$), hence for fixed $m$ and $0<\phi_2<\infty$,

\[
  \sqrt{N}\,\Bigl\{Q^{(m)}_j
                   -Q^{(m),\circ}_j\Bigr\}
  \ \xrightarrow{p}\ 0, \qquad j=1,2.
\]
Therefore we may work with $Q^{(m),\circ}_1,Q^{(m),\circ}_2$, which are both quadratic forms
of the same field $D^{(m)}_{[1]}X$ over the same index set, without affecting the
$\sqrt N$ limits. Note that $\E Q^{(m)}_j = \E Q^{(m),\circ}_j$ by stationarity.

\subsection{ID CLTs for \texorpdfstring{$(Q^{(m)}_{1},Q^{(m)}_{2})$} and the estimators}
\label{subsec:IDCLTm}

\begin{proposition}[Quadratic–form CLT, higher order]\label{thm:UniCLTm}
Let $\phi_{2}\in(k,k+1)$ and $m\ge k+1$.  For every $(\alpha,\beta)\neq(0,0)$,
\[
   \sqrt{N}\,\Bigl\{
     \alpha\bigl(Q^{(m)}_{1}-q^{(m)}_{1}\bigr)
    +\beta \bigl(Q^{(m)}_{2}-q^{(m)}_{2}\bigr)\Bigr\}
   \xrightarrow{d}
   \mathcal N\!\bigl(0,\sigma^{2}_{\alpha,\beta;m}\bigr),
\]
with
\[
   \sigma^{2}_{\alpha,\beta;m}
   = 2\int_{\T^{2}}
        \bigl(\alpha+\beta B_{m}(\lambda)\bigr)^{2}
        F(\lambda)^{2}\, \mu(\mathrm d\lambda).
\]
\end{proposition}
\begin{proof}
Refer to \ref{proofs}
\end{proof}

\begin{theorem}[Joint CLT and delta–method]\label{cor:Jointm}
Let $\Sigma^{(m)}$ have entries
\[
 \Sigma^{(m)}_{\ell r}
   = 2\int_{\T^{2}}
           b_{\ell}(\lambda)b_{r}(\lambda)\,
           F(\lambda)^{2}\,\mu(\mathrm d\lambda),
 \qquad
 b_{1}\equiv 1,\quad b_{2}=B_{m}.
\]
Then
\[
   \sqrt{N}
   \begin{pmatrix}
     Q^{(m)}_{1}-q^{(m)}_{1}\\[2pt]
     Q^{(m)}_{2}-q^{(m)}_{2}
   \end{pmatrix}
   \xrightarrow{d}\mathcal N(\bm0,\Sigma^{(m)}).
\]
Let
\(
   h^{(m)}(y_{1},y_{2})
   =\bigl(y_{1}/a_{m}(\tfrac12\log_{2}(y_{2}/y_{1})),
          \tfrac12\log_{2}(y_{2}/y_{1})\bigr)
\)
and $J^{(m)}=\nabla h^{(m)}(q^{(m)}_{1},q^{(m)}_{2})$. Then
\[
   \sqrt{N}
   \begin{pmatrix}
     \widehat\phi_{1}-\phi_{1}\\[2pt]
     \widehat\phi_{2}-\phi_{2}
   \end{pmatrix}
   \xrightarrow{d}
   \mathcal N\!\bigl(\bm0,\,J^{(m)}\Sigma^{(m)}J^{(m)\top}\bigr).
\]
\end{theorem}
\begin{proof}
Refer to \ref{proofs}
\end{proof}

\begin{remark}[Normalisation by $N$]\label{rem:NormByN}
Since $|\Lambda_{n-2jm}|/N\to1$ as $n\to\infty$, replacing $|\Lambda_{n-2jm}|$ by $N$
throughout leaves all $\sqrt{N}$ limits and asymptotic covariances unchanged (when centering
at their own means), exactly as in the IRF–0 case.
\end{remark}

\subsection{FD CLTs via exact rescaling for m-th differences}\label{subsec:FDm}

Fix $m\ge \lfloor\phi_2\rfloor+1$ and let $N=n^2$. In the FD design on $[0,1]^2$,
define the FD parameters
\[
  \tau_1 := \phi_1\,N^{-\phi_2},\qquad \tau_2 := \phi_2,\qquad \tau=(\tau_1,\tau_2).
\]
By the spatial reindexing and amplitude rescaling argument in
Section~\ref{sec:IDtoFD} (see \eqref{eq:fd-id-ratio-1}),
and since $D^{(m)}$ is linear while each $Q^{(m)}_j$ is homogeneous of degree two,
we have the exact distributional identity (for the same $m$):
\begin{equation}\label{eq:fd-id-ratio-m}
  \left(\frac{\widehat\tau_{1}^{(m)}}{\tau_1},\ \widehat\tau_{2}^{(m)}\right)
  \stackrel{d}{=} 
  \left(\frac{\widehat\phi_{1}^{(m)}}{\phi_1},\ \widehat\phi_{2}^{(m)}\right).
\end{equation}

\paragraph{ID ratio map and Jacobian (order $m$).}
Let
\[
\bm Q^{(m)}:=
\begin{pmatrix}Q^{(m)}_{1}\\ Q^{(m)}_{2}\end{pmatrix},
\qquad
\bm q^{(m)}:=
\begin{pmatrix}q^{(m)}_{1}\\ q^{(m)}_{2}\end{pmatrix},
\]
and denote by \(\Sigma^{(m)}\) the covariance from the joint quadratic–form CLT
for order \(m\) (Corollary~\ref{cor:Jointm}):
\[
\sqrt{N}\,\big(\bm Q^{(m)}-\bm q^{(m)}\big)\Rightarrow \mathcal N\!\big(\bm 0,\Sigma^{(m)}\big).
\]
Define the delta map for order \(m\),
\[
h_m(y_1,y_2)
=\left(
\frac{y_1}{a_m\!\big(\tfrac12\log_2(y_2/y_1)\big)},
\ \tfrac12\log_2(y_2/y_1)
\right),
\]
and its Jacobian at \(\bm q^{(m)}\),
\[
J^{(m)}:=\nabla h_m\big(q^{(m)}_{1},q^{(m)}_{2}\big).
\]
Writing \(A_m(\phi):=a_m(\phi)\), \(s:=(2\ln 2)^{-1}\), and
\(r:=q^{(m)}_{1}/q^{(m)}_{2}=2^{-2\phi_2}\) (the ratio does not depend on \(m\)),
we have (as in the first–order case with \(a_1\) replaced by \(a_m\))
\[
J^{(m)}=
\begin{pmatrix}
A_m^{-1}+sA_m'/A_m^2 & -\,sA_m'/A_m^2\,r\\[4pt]
-\,s/q^{(m)}_{1} & \ \ s/q^{(m)}_{2}
\end{pmatrix},
\qquad
q^{(m)}_{1}=\phi_1 A_m(\phi_2),\quad q^{(m)}_{2}=2^{2\phi_2}q^{(m)}_{1}.
\]

To obtain the result for \(\big(\widehat\phi_1/\phi_1-1,\ \widehat\phi_2-\phi_2\big)\),
compose with \(g(u,v)=(u/\phi_1-1,\ v-\phi_2)\) and set
\[
G_m:=g\circ h_m,\qquad
K^{(m)}:=\nabla G_m(\bm q^{(m)})=\begin{pmatrix}1/\phi_1&0\\[2pt]0&1\end{pmatrix} J^{(m)}.
\]
Hence
\[
K^{(m)}=
\begin{pmatrix}
\dfrac{1}{\phi_1}\!\left(A_m^{-1}+sA_m'/A_m^2\right) & -\,\dfrac{s}{\phi_1}\,\dfrac{A_m'}{A_m^2}\,r\\[10pt]
-\,\dfrac{s}{q^{(m)}_{1}} & \ \ \dfrac{s}{q^{(m)}_{2}}
\end{pmatrix}.
\]

\paragraph{ID ratio CLT (order $m$).}
By the multivariate delta method applied to Corollary~\ref{cor:Jointm},
\[
\sqrt{N}
\begin{pmatrix}
\widehat\phi_{1}^{(m)}/\phi_1-1\\[2pt]
\widehat\phi_{2}^{(m)}-\phi_2
\end{pmatrix}
\ \Rightarrow\ 
\mathcal N\!\big(\bm 0,\ K^{(m)}\,\Sigma^{(m)}\,{K^{(m)}}^\top\big).
\]

\paragraph{FD ratio CLT via exact rescaling.}
From \eqref{eq:fd-id-ratio-m} we obtain the exact FD counterpart:
\begin{equation}\label{eq:FDm-ratio-CLT}
\sqrt{N}
\begin{pmatrix}
\widehat\tau_{1}^{(m)}/\tau_1-1\\[6pt]
\widehat\tau_{2}^{(m)}-\tau_2
\end{pmatrix}
\ \Rightarrow\ 
\mathcal N\!\big(\bm 0,\ K^{(m)}\,\Sigma^{(m)}\,{K^{(m)}}^\top\big).
\end{equation}
\emph{Equivalently} (expressed in terms of $(\phi_1,\phi_2)$),
\begin{equation}\label{eq:FDm-equivalent-phi}
\sqrt{N}
\begin{pmatrix}
\dfrac{N^{\widehat\phi_{2}^{(m)}}\,\widehat\tau_{1}^{(m)}}{\phi_1}-1\\[8pt]
\widehat\tau_{2}^{(m)}-\phi_2
\end{pmatrix}
\ \Rightarrow\ 
\mathcal N\!\big(\bm 0,\ K^{(m)}\,\Sigma^{(m)}\,{K^{(m)}}^\top\big).
\end{equation}

\begin{figure}[!htbp]
  \centering
  \includegraphics[width=0.8\linewidth]{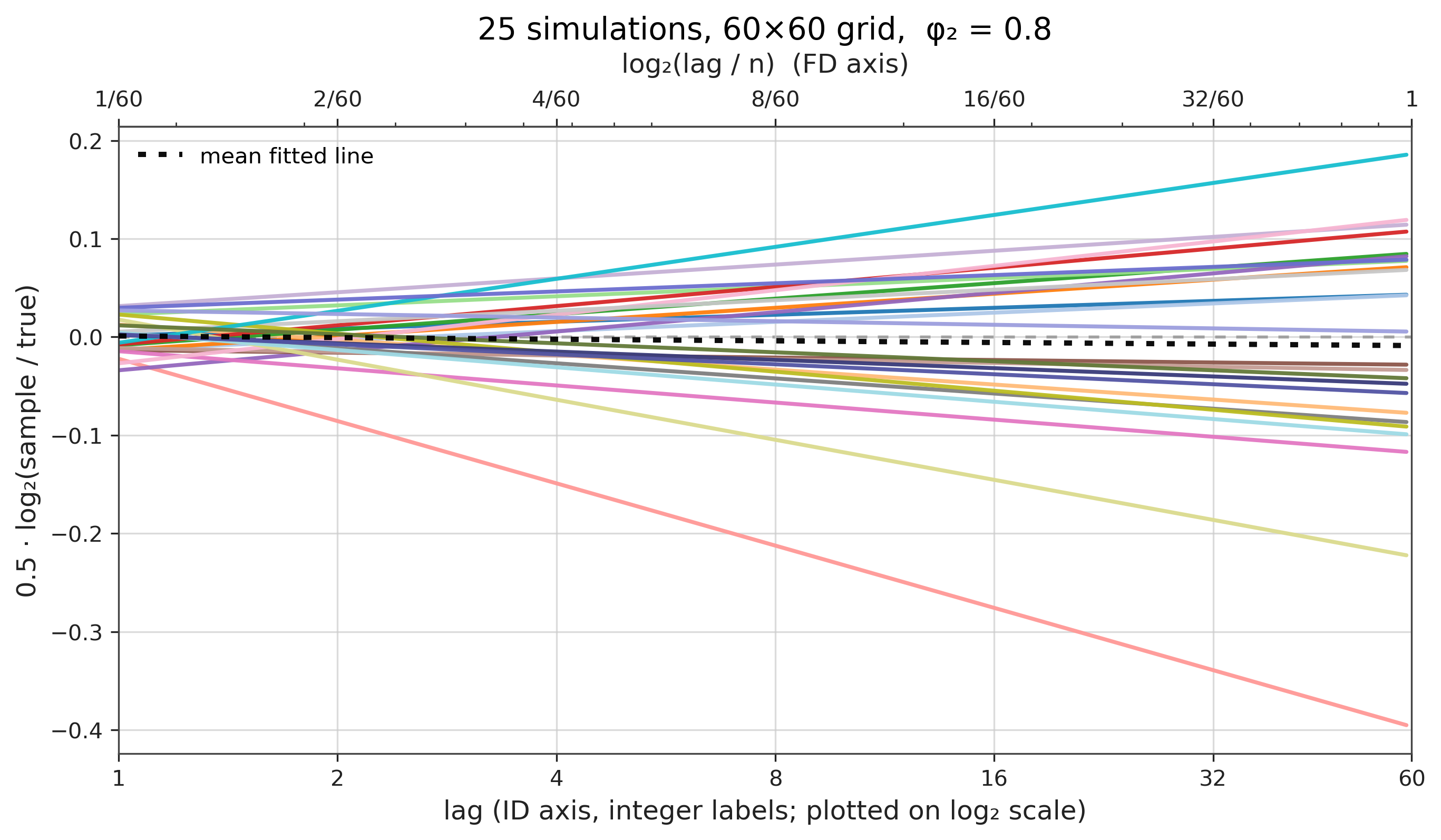}%
  \caption{\textbf{Sample vs.\ true variograms on a log scale (bilinear, IRF-0).} For 25 simulated \(60\times60\) power-law fields with \(\phi_2=0.8\), the plot compares estimated and true variograms on a log--log scale. Each colored line shows the bilinear ratio slope \(\hat{\phi}_2-\phi_2\); a mean slope near zero indicates unbiasedness and the spread reflects sampling variability.}
  \label{fig:straight-variogram}
\end{figure}

Figure \ref{fig:straight-variogram} makes the core invariance behind our method visible and shows why fixed- and increasing-domain designs are effectively the same once distance is measured in ``neighbor units.'' 
For each simulation, we plot the discrepancy 
$\tfrac12\log_{2}\!\big\{\hat\gamma(h)/\gamma_{\text{true}}(h)\big\}$ 
against $\log_{2}h$. For a power law $\gamma(h)\propto h^{2\phi_{2}}$, this curve is a straight line with slope $\hat\phi_{2}-\phi_{2}$.  Because the bilinear two–scale estimator uses a fixed lag ratio $2{:}1$, $\hat\phi_{2}$ is determined entirely by that slope and is therefore unit–free; $\hat\phi_{1}$ then follows by matching the level at $h=1$. The bottom axis labels integer lags $h$ (ID view), while the top axis relabels the same points as $h/n$ (FD view); this is a horizontal shift by $-\log_{2}n$ that leaves slopes unchanged. Thus rescaling the grid only re–anchors the \emph{level} (i.e., $\phi_{1}$), while the \emph{roughness} $\phi_{2}$ is invariant precisely the equivalence used to transfer ID CLTs to FD after the exact rescaling. The dashed horizontal line marks perfect agreement $(\hat\gamma=\gamma_{\text{true}})$ and the dotted black line is the mean fitted line across simulations; its near–zero slope indicates approximate unbiasedness for $\phi_{2}$ here, and the vertical spread visualizes sampling variability. Figure \ref{fig:straight-variogram} makes the FD asymptotic results for $(\widehat\phi_1,\widehat\phi_2)$ qualitatively clear. Let us focus on the distance scale on the bottom axis, so that the distance between neighboring observations is 1.
At this lag, $\widehat\gamma(1) = \widehat\phi_1-\Gamma(-\widehat\phi_2)$, which, from the ID asymptotics, has relative variability of order $N^{-1/2}$. 
For lag $n$ at the right edge of the plot, we have 
\begin{equation}
    \log\left(\frac{\widehat\gamma(n)}{\gamma_{\mathrm{true}}(n)}\right) = \log\left(\frac{\widehat\phi_1}{\phi_1}\right) + (\widehat\phi_2-\phi_2)\log n-\log\left(\frac{\Gamma(-\widehat\phi_2)}{\Gamma(-\phi_2)}\right).
    \label{eq:gamma_n}
\end{equation}
Because both $\widehat\phi_1$ and $\widehat\phi_2$ have errors of order $N^{-1/2}$, the presence of the $\log n$ factor multiplying $\widehat\phi_2-\phi_2$ in (\ref{eq:gamma_n}) immediately implies that $\log(\widehat\gamma(n)/\gamma_{\mathrm{true}}(n))$ is of order $N^{-1/2}\log N$ in probability and, furthermore, $\log\widehat\gamma(n)$ and $\widehat\phi_2$ are strongly correlated.
But, since using the FD distance scale in the upper axis,
$\log(\widehat\phi_1/\phi_1) = \log\widehat\gamma(n) + \log(\Gamma(-\widehat\phi_2)/\Gamma(-\phi_2))$,
the $N^{-1/2}\log N$ convergence rate for $\log(\widehat\phi_1/\phi_1)$ and the strong correlation between $\widehat\phi_1$ and $\widehat\phi_2$ under FD asymptotics follows.

\begin{figure*}[!htbp]
  \centering
  \includegraphics[width=0.75\textwidth]{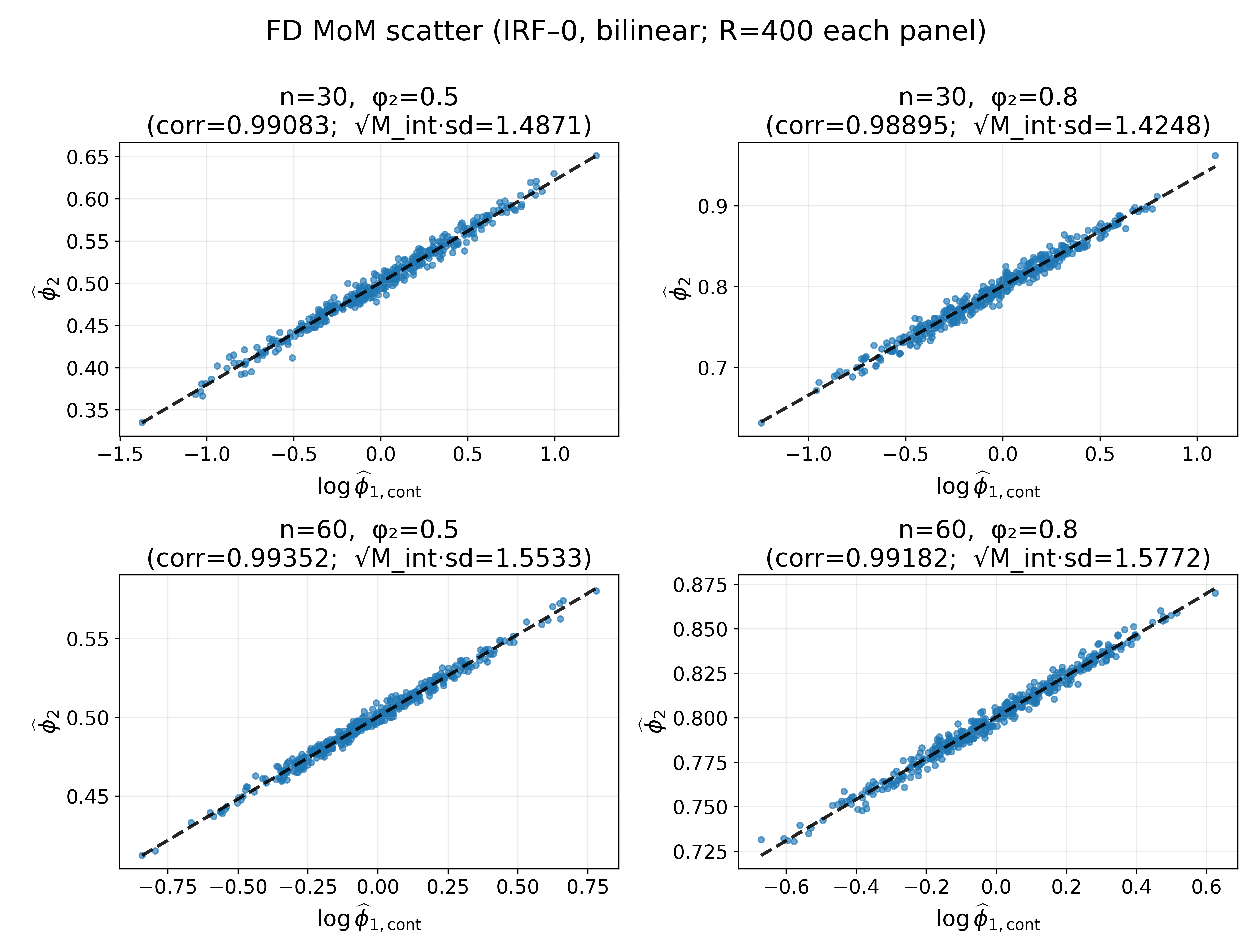}

  \caption{\textbf{Fixed–domain (FD) estimator scatter plot for IRF–0 (bilinear).} Each panel shows 400 replicates of \((\log \widehat{\phi}_1,\widehat{\phi}_2)\) from the bilinear method-of-moments on an \(n\times n\) grid. The near-linear ridge (\(\mathrm{corr}\approx0.99\)) and stable scaled s.d.\ (\(\approx1.45\)) indicate \(\widehat{\phi}_2\) is unbiased and that \(\operatorname{sd}(\widehat{\phi}_2)\approx c/\sqrt{M_{\mathrm{int}}}\) with \(c\approx1.5\).}
  \label{fig:fd-scatter}
\end{figure*}

\section{Robustness to irregular sampling}
\label{sec:MatrixForm}

In practice, lattice observations are often incomplete or slightly perturbed.
We study the stability of the quadratic–variation statistics under four departures from the ideal design:
(i) \emph{deterministic deletions} of a vanishingly small fraction of sites,
(ii) \emph{Bernoulli thinning} with high retention probability,
(iii) \emph{small spatial jitter} in the fixed–domain (FD) setting, and
(iv) \emph{model misspecification}: Matérn truth with a power–law working model.
In all cases, the perturbations alter the quadratic forms only at lower order, so the CLTs and the method–of–moments inference for $(\phi_1,\phi_2)$ remain valid with the same asymptotic covariance matrix.

Throughout fix an integer \(m\ge \lfloor\phi_2\rfloor+1\), and use the two stencils
\(D^{(m)}_{[1]}\) and \(D^{(m)}_{[2]}\) from \S\ref{subsec:DiffM}.

\subsection{Quadratic forms on the lattice and uniform bounds}
\label{subsec:FilterMatrices}

Let $\Lambda_n=\{0,\dots,n-1\}^2$ with $N=n^2$ and $M_j=(n-jm)^2$ for $j=1,2$.
Write $\mathbf X_n=(X_t)_{t\in\Lambda_n}\in\R^{N}$.
Let $F^{(m)}_{[j]}$ be the $(M_j\times N)$ matrix that stacks \(D^{(m)}_{[j]}X_t\) over $t\in\Lambda_{n-jm}$, so
\[
  \bigl(D^{(m)}_{[j]}X\bigr) = F^{(m)}_{[j]}\,\mathbf X_n,\qquad
  Q^{(m)}_{j}=\frac{1}{M_j}\|F^{(m)}_{[j]}\mathbf X_n\|_2^2
             = \mathbf X_n^\top A^{(m)}_{j}\mathbf X_n,
\]
with $A^{(m)}_{j}=M_j^{-1}(F^{(m)}_{[j]})^\top F^{(m)}_{[j]}$ and $\|A^{(m)}_{j}\|_2\le C_m$ (constant depends only on $m$).
As in \S\ref{subsec:TwoScaleM}, $F^{(m)}_{[2]}=H^{(m)}F^{(m)}_{[1]}$, where $H^{(m)}$ is the $(m{+}1)^2$–point block–sum matrix with $\|H^{(m)}\|_2\le 4^m$; hence
\[
  A^{(m)}_{2}=\frac{1}{(n-2m)^2}\,(F^{(m)}_{[1]})^\top (H^{(m)})^\top H^{(m)}F^{(m)}_{[1]}.
\]

Let
\[
  \mathbf D^{(m)}_{[1],n} := \bigl(D^{(m)}_{[1]}X_t\bigr)_{t\in\Lambda_{n-m}}\in\R^{M_1},\qquad
  \Sigma^{(m)}_{n}=\Cov\!\bigl(\mathbf D^{(m)}_{[1],n}\bigr).
\]
By Lemma~\ref{lem:L2m}, $f_{D^{(m)}_{[1]}X}(\lambda)=|g_m(\lambda)|^2 f_0^{\lat}(\lambda)$ is bounded on $\T^2$; thus
\begin{equation}\label{eq:SigmaBound}
  \sup_{n}\|\Sigma^{(m)}_{n}\|_2 \ \le\ \bar C_m<\infty,
  \qquad \|\Sigma^{(m)}_{n}\|_F\ \lesssim\ \sqrt N.
\end{equation}

For $(\alpha,\beta)\in\R^2$ define the bounded kernel on the unit–step outputs
\begin{equation}\label{eq:BalphaBeta}
  B_{\alpha,\beta}
  := \frac{\alpha}{M_1}\,I_{M_1}
     + \frac{\beta}{M_2}\,(H^{(m)})^\top H^{(m)}.
\end{equation}
Since $M_j\asymp n^2$ and $\|H^{(m)}\|_2\le 4^m$, we have uniform bounds
\begin{equation}\label{eq:Bnorms}
  \|B_{\alpha,\beta}\|_2 \ \lesssim\ n^{-2}(|\alpha|+|\beta|),
  \qquad
  \|B_{\alpha,\beta}\|_F \ \lesssim\ n^{-2}\sqrt N\,(|\alpha|+|\beta|).
\end{equation}
Finally note
\begin{equation}\label{eq:QasQuadForm}
  \alpha Q^{(m)}_1+\beta Q^{(m)}_2
  \;=\; \mathbf D^{(m)\top}_{[1],n}\,B_{\alpha,\beta}\,\mathbf D^{(m)}_{[1],n}.
\end{equation}

\subsection{Perturbation lemmas}
\label{subsec:MasterPerturb}

We collect two statistical perturbation results for quadratic forms that we
reuse across deletions, thinning, jitter, and misspecification.
Let $Y_n:=\mathbf D^{(m)}_{[1],n}$ denote the filtered data vector
and $B_{\alpha,\beta}$ the (deterministic) weight matrix for the quadratic form,
as defined in \eqref{eq:BalphaBeta}. We also use the uniform bounds
\eqref{eq:SigmaBound}–\eqref{eq:Bnorms}.

For a centered Gaussian vector $Y$ with covariance $\Sigma$ and a symmetric matrix $B$,
\[
  \E(Y^\top BY)=\mathrm{tr}(\Sigma B),\qquad
  \Var(Y^\top BY)=2\,\mathrm{tr}\bigl((\Sigma B)^2\bigr).
\]
Here $\| \cdot\|_F$ is Frobenius norm.

\begin{lemma}[Perturbation of weight and covariance]\label{lem:MasterKernel}
Let $B_n:=B_{\alpha,\beta}$ and $\Sigma_n:=\Sigma^{(m)}_n$.
Consider perturbed objects $B'_n=B_n+\Delta B_n$ and $\Sigma'_n=\Sigma_n+\Delta\Sigma_n$, with $\Delta B_n,\Delta\Sigma_n$ symmetric. With constants depending only on $m$ and $(\alpha,\beta)$:

\emph{(a) Weight (kernel) perturbation.} Expectations and variances computed under $\Sigma_n$ satisfy
\[
\bigl|\E(Y_n^\top B'_n Y_n)-\E(Y_n^\top B_n Y_n)\bigr|
 \ \le\ \bar C_m\,\|\Delta B_n\|_F,
\]
\[
\bigl|\Var(Y_n^\top B'_n Y_n)-\Var(Y_n^\top B_n Y_n)\bigr|
 \ \lesssim\ \|B_n\|_F\,\|\Delta B_n\|_F+\|\Delta B_n\|_F^2.
\]

\emph{(b) Covariance perturbation.} Writing $\E'$ and $\Var'$ for expectation/variance under $\Sigma'_n$,
\[
\bigl|\E(Y_n^\top B_n Y_n)-\E'(Y_n^\top B_n Y_n)\bigr|
 \ \le\ \|B_n\|_F\,\|\Delta\Sigma_n\|_F,
\]

\[
\bigl|\Var(Y_n^\top B_n Y_n)-\Var'(Y_n^\top B_n Y_n)\bigr|
 \ \lesssim\ \|B_n\|_2^2\Bigl(\|\Delta\Sigma_n\|_F\,\|\Sigma_n\|_F+\|\Delta\Sigma_n\|_F^2\Bigr).
\]

In particular, by \eqref{eq:SigmaBound}–\eqref{eq:Bnorms},
\[
 \|B_n\|_F\ \lesssim\ n^{-2}\sqrt N,\qquad
 \|B_n\|_2\ \lesssim\ n^{-2},\qquad
 \|\Sigma_n\|_F\ \lesssim\ \sqrt N.
\]
\end{lemma}

\begin{proof}
Use the Gaussian identities above and the trace/Hölder inequalities
$\bigl|\mathrm{tr}(A^\top C)\bigr|\le \|A\|_F\|C\|_F$ and
$\|XY\|_F\le \|X\|_2\|Y\|_F$, then insert \eqref{eq:SigmaBound}–\eqref{eq:Bnorms}.
\end{proof}

\begin{lemma}[Perturbation by observation error]\label{lem:MasterInput}
Let $Y_n:=\mathbf D^{(m)}_{[1],n}$ and suppose we observe $Y_n^\star=Y_n+\delta_n$,
with $\max_{t}\E\,\delta_{n,t}^2 \le u_n^2\downarrow 0$. For any fixed $(\alpha,\beta)$,
\[
\E\bigl|\,Y_n^{\star\top}B_{\alpha,\beta}Y_n^\star - Y_n^\top B_{\alpha,\beta}Y_n\,\bigr|
\ \lesssim\ u_n,
\]
and
\[
Y_n^{\star\top}B_{\alpha,\beta}Y_n^\star - Y_n^\top B_{\alpha,\beta}Y_n
\ =\ O_p\!\bigl(n^{-2}\sqrt N\,u_n\bigr)\ \ \vee\ \ O_p(u_n^2).
\]
\end{lemma}

\begin{proof}
Expand $2Y_n^\top B_{\alpha,\beta}\delta_n+\delta_n^\top B_{\alpha,\beta}\delta_n$ and apply Cauchy–Schwarz with \eqref{eq:Bnorms}–\eqref{eq:SigmaBound}.
\end{proof}

\subsection{Deterministic sparse deletions}\label{sec:SparseDelete}

Let $\mathcal M_n\subset\Lambda_n$ be a set of deleted sites of size $k_n$.
Remove from $F^{(m)}_{[j]}$ every row that touches a deleted site to obtain
$\widetilde F^{(m)}_{[j]}$ with $\widetilde M_j$ rows (so $M_j-\widetilde M_j\le c_{j,m}k_n$).
Define the pruned quadratic forms with their own divisors:
\[
  \widetilde Q^{(m)}_{j}
  =\frac{1}{\widetilde M_j}\bigl\|\widetilde F^{(m)}_{[j]}\mathbf X_n\bigr\|_2^2
  = \mathbf X_n^\top \widetilde A^{(m)}_{j}\mathbf X_n,\qquad
  \widetilde A^{(m)}_{j}=\frac{1}{\widetilde M_j}\,
      (\widetilde F^{(m)}_{[j]})^\top \widetilde F^{(m)}_{[j]} .
\]

\begin{lemma}[Matrix perturbation by row deletions]\label{lem:opnorm}
There exists $C_m<\infty$ (independent of $n$) such that, for $j=1,2$,
\[
  \bigl\|\widetilde A^{(m)}_{j}-A^{(m)}_{j}\bigr\|_2
  + \bigl\|\widetilde A^{(m)}_{j}-A^{(m)}_{j}\bigr\|_F
  \ \le\ C_m\,\frac{k_n}{N},
  \qquad
  \mathrm{rank}\bigl(\widetilde A^{(m)}_{j}-A^{(m)}_{j}\bigr)\ \le\ C_m\,k_n .
\]
Consequently,
\(
  \|\widetilde A^{(m)}_{\alpha,\beta}-A^{(m)}_{\alpha,\beta}\|_{2,F}
 \le C_m(|\alpha|+|\beta|)\,k_n/N.
\)
\end{lemma}

\begin{proof}
As in the $m{=}1$ case, using $\|F^{(m)}_{[j]}\|_F^2=C_m M_j$ and that each removed row $r_\ell r_\ell^\top$ has $\|r_\ell\|_2^2=C_m$.
\end{proof}

Let $S_j$ be the row–selection matrix keeping rows not touching deleted sites, so
$\widetilde F^{(m)}_{[1]}=S_1 F^{(m)}_{[1]}$ and $\widetilde F^{(m)}_{[2]}=S_2 H^{(m)} F^{(m)}_{[1]}$.
Then, on the space of $Y_n=\mathbf D^{(m)}_{[1],n}$,
\[
  \widetilde B_{\alpha,\beta}
  := \frac{\alpha}{\widetilde M_1}\,S_1^\top S_1
   + \frac{\beta}{\widetilde M_2}\,(H^{(m)})^\top S_2^\top S_2 H^{(m)},
\quad
  \Delta B_{\alpha,\beta}:=\widetilde B_{\alpha,\beta}-B_{\alpha,\beta},
\]
satisfy $\|\Delta B_{\alpha,\beta}\|_{2,F}=O(k_n/N)$ by Lemma~\ref{lem:opnorm} and $\|H^{(m)}\|_2\le 4^m$.

\begin{lemma}[Stability of quadratic forms]\label{lem:SparsePerturb}
Under $\sup_n\|\Sigma^{(m)}_{n}\|_2<\infty$, for any fixed $(\alpha,\beta)$,
\[
  \E\Bigl[\mathbf X_n^\top(\widetilde A^{(m)}_{\alpha,\beta}-A^{(m)}_{\alpha,\beta})\mathbf X_n\Bigr]
  = O\!\left(\frac{k_n}{N}\right),\qquad
  \Var\Bigl[\mathbf X_n^\top(\widetilde A^{(m)}_{\alpha,\beta}-A^{(m)}_{\alpha,\beta})\mathbf X_n\Bigr]
  = O\!\left(\frac{k_n^2}{N^2}\right).
\]
Consequently,
\(
  \mathbf X_n^\top(\widetilde A^{(m)}_{\alpha,\beta}-A^{(m)}_{\alpha,\beta})\mathbf X_n
  = O_p(k_n/N).
\)
\end{lemma}

\begin{proof}
Apply Lemma~\ref{lem:MasterKernel} with $\Delta B_n=\Delta B_{\alpha,\beta}$ and \eqref{eq:Bnorms}–\eqref{eq:SigmaBound}.
\end{proof}

\begin{corollary}[CLTs with deterministic gaps]\label{cor:SparseCLT}
If $k_n=o(\sqrt N)$, then for any fixed $(\alpha,\beta)$,
\[
  \sqrt{N}\Bigl\{\mathbf X_n^\top\widetilde A^{(m)}_{\alpha,\beta}\mathbf X_n
                     -\E[\mathbf X_n^\top\widetilde A^{(m)}_{\alpha,\beta}\mathbf X_n]\Bigr\}
  - \sqrt{N}\Bigl\{\mathbf X_n^\top A^{(m)}_{\alpha,\beta}\mathbf X_n
                     -\E[\mathbf X_n^\top A^{(m)}_{\alpha,\beta}\mathbf X_n]\Bigr\}
  \xrightarrow{p} 0 .
\]
Hence the (ID and FD) joint CLTs for $(Q^{(m)}_{1},Q^{(m)}_{2})$ continue to hold for
$(\widetilde Q^{(m)}_{1},\widetilde Q^{(m)}_{2})$ with the same asymptotic covariance.
\end{corollary}

\subsection{Bernoulli thinning}\label{sec:Bernoulli}

Suppose each site is retained independently with probability $p_n\in(0,1]$.
Let $\mathcal M_n$ be the deleted set; then $\E k_n=N(1-p_n)$ and $\Var(k_n)=Np_n(1-p_n)$.

\begin{assumption}[High retention]\label{ass:pn-correct}
\[
 (1-p_n)\sqrt{N}\ \longrightarrow\ 0\qquad\text{(ID and FD)}.
\]
\end{assumption}

\begin{lemma}\label{lem:mn-small}
Under Assumption~\ref{ass:pn-correct}, $k_n=o_p(\sqrt N)$.
\end{lemma}

\begin{proof}
Markov’s inequality: $\Pr(k_n>\varepsilon \sqrt N)\le (1-p_n)\sqrt N/\varepsilon\to 0$.
\end{proof}

\begin{theorem}[CLTs under Bernoulli thinning]\label{thm:Bernoulli}
Under Assumption~\ref{ass:pn-correct}, the ID and FD CLTs and the delta–method limits for $(\widehat\phi_1,\widehat\phi_2)$ remain valid with the same asymptotic covariances.
\end{theorem}

\begin{proof}
Condition on the thinning pattern and apply Lemmas~\ref{lem:opnorm}, \ref{lem:mn-small}, and \ref{lem:MasterKernel}, then Slutsky.
\end{proof}

\subsection{Small spatial jitter (FD)}\label{sec:Jitter}

Let the FD sampling locations be \(x_t=t/n+\varepsilon_t\) with \(\|\varepsilon_t\|_\infty\le c/n\) (\(c>0\)).
We observe
\(
  X_t := X(t/n),\ X_t^\star := X(x_t).
\)
For \(j\in\{1,2\}\) define the jitter–induced perturbations
\[
   \Delta^{(m)}_{j,t}:= D^{(m)}_{[j]} X^\star_t - D^{(m)}_{[j]} X_t,\qquad
   \Delta Q^{(m)}_{j}:= Q^{(m)}_{j}(X^\star)-Q^{(m)}_{j}(X).
\]

\begin{lemma}[Second–moment control]\label{lem:jitter-L2}
There exists $C_m<\infty$ such that
\(
  \sup_{t} \E\,\bigl(\Delta^{(m)}_{j,t}\bigr)^2 \le C_m\, n^{-2\phi_2}
\)
for $j=1,2$.
\end{lemma}

\begin{proof}
Each $\Delta^{(m)}_{j,t}$ is a finite linear combination of increments $X(x)-X(y)$ with $\|x-y\|\lesssim n^{-1}$. Use the semivariogram bound $2\gamma_{\phi_1,\phi_2}(h)\lesssim \|h\|^{2\phi_2}$.
\end{proof}

\begin{proposition}[Effect on $Q^{(m)}_{j}$]\label{prop:jitter-Q}
For $j=1,2$,
\[
   \E\bigl|\Delta Q^{(m)}_{j}\bigr|\ \lesssim\ n^{-\phi_2},
   \qquad
   \sqrt{N}\,\Delta Q^{(m)}_{j} \xrightarrow{p}0 \quad\text{if }\ \phi_2>1.
\]
\end{proposition}

\begin{proof}
Because $D^{(m)}$ is linear, the difference lives at the output level:
$Y_n^\star=Y_n+\delta_n$ with $Y_n=\mathbf D^{(m)}_{[1],n}$ and $\delta_n=(\Delta^{(m)}_{1,t})_t$.
By Lemma~\ref{lem:jitter-L2}, $u_n\asymp n^{-\phi_2}$; apply Lemma~\ref{lem:MasterInput}.
This gives $\E|\Delta Q^{(m)}_j|\lesssim u_n$ and
$\Delta Q^{(m)}_j=O_p(n^{-2}\sqrt N\,u_n)\vee O_p(u_n^2)=O_p(n^{-\phi_2})$.
Hence $\sqrt N\,\Delta Q^{(m)}_{j}\to 0$ when $\phi_2>1$.
\end{proof}

\begin{corollary}[FD CLTs with jitter]\label{cor:jitter-CLT}
If $\phi_2>1$ and $\|\varepsilon_t\|_\infty\le c/n$, the FD $\sqrt{N}$–CLTs for
\(
  \bigl(n^{2\widehat\phi_{2,n}^{\mathrm{FD}}}\widehat\phi_{1,n}^{\mathrm{FD}},
        \widehat\phi_{2,n}^{\mathrm{FD}}\bigr)
\)
remain valid with the same asymptotic covariances as in ~\ref{cor:Jointm}.
\end{corollary}

\subsection{Robustness under fixed–domain misspecification (Matérn truth)}
\label{sec:Misspecification}

We study the FD behaviour of the two–scale, order-$m$ estimators when the
\emph{true} field is Matérn but we fit the power–law (PL) working equations.
Work on $[0,1]^2$ with observations $\{X(t/n):t\in\Lambda_n\}$, $N=n^2$.
Keep $D^{(m)}_{[j]}$, $Q^{(m)}_j$, $F^{(m)}_{[j]}$, $A^{(m)}_j$
from \S\ref{subsec:DiffM}–\S\ref{subsec:FilterMatrices}.

\paragraph{Matérn covariance and spectrum.}
For $\sigma^2>0$, $\nu\in(0,1)$ and $\kappa=\sqrt{2\nu}/\rho$,
\[
 C_{\Mat}(h)=\sigma^2\frac{2^{1-\nu}}{\Gamma(\nu)}(\kappa\|h\|)^\nu K_\nu(\kappa\|h\|),
\]
and, with the convention $k(h)=\int_{\R^d}e^{2\pi i s\cdot h}S(s)\,ds$,
\[
 S_{\Mat}(s)=\sigma^2\,\frac{\Gamma(\nu+\tfrac d2)}{\Gamma(\nu)\,\pi^{d/2}}\,
              \kappa^{2\nu}\bigl(\kappa^2+4\pi^2\|s\|^2\bigr)^{-(\nu+d/2)}.
\]

\paragraph{Small–lag expansion and tangent PL.}
In $d=2$ and for $0<\nu<1$,
\begin{equation}\label{eq:Mat-small-lag}
  \gamma_{\Mat}(r)= c_{\Mat}\,r^{2\nu} + a_{\Mat}\,r^{2} + O(r^{2\nu+2}+r^{4}),
\quad
c_{\Mat}=\sigma^2\,\frac{\kappa^{2\nu}}{2^{2\nu}}\,\frac{|\Gamma(-\nu)|}{\Gamma(\nu)},\ \
a_{\Mat}=\sigma^2\,\frac{\kappa^2}{4(1-\nu)}.
\end{equation}
The \emph{tangent} PL at the origin is therefore the IRF-0 PL model with
$(\phi_1,\phi_2)=(c_{\Mat},\nu)$; see \cite[Sec.~2.10–2.11]{Stein1999}.

\paragraph{Filtered spectra on the torus.}
Let $f_{\PL,\nu}^{\lat}$ denote the aliased PL (tangent) spectrum and put
\[
\widetilde f^{(m)}_{\Mat}(\lambda)=|g_m(\lambda)|^2 f^{\lat}_{\Mat}(\lambda),
\qquad
\widetilde f^{(m)}_{\PL,\nu}(\lambda)=|g_m(\lambda)|^2 f^{\lat}_{\PL,\nu}(\lambda),
\qquad \lambda\in\T^2,
\]
where $g_m$ is the $m$th-order differencing filter and $f^{\lat}$ is the usual periodization.

\begin{lemma}[Filtered spectral domination]\label{lem:domination}
For $d=2$, $0<\nu<1$ and any $m\ge1$,
\[
  \widetilde f^{(m)}_{\Mat}(\lambda)\ \le\ \widetilde f^{(m)}_{\PL,\nu}(\lambda)
  \quad\text{for all }\lambda\in\T^2,
\]
and $r_m(\lambda):=\widetilde f^{(m)}_{\PL,\nu}(\lambda)-\widetilde f^{(m)}_{\Mat}(\lambda)$
belongs to $L^1(\T^2)\cap L^2(\T^2)$ whenever $m-\nu>\tfrac12$.
\end{lemma}

\begin{proof}
For continuum frequencies $s\in\R^2$ we have
\(
(\kappa^2+4\pi^2\|s\|^2)^{-(\nu+1)}
\le (4\pi^2\|s\|^2)^{-(\nu+1)}
\)
so $S_{\Mat}(s)\le C_\nu\,\|s\|^{-2\nu-2}$ with the \emph{same} leading constant
as the tangent PL. Periodization preserves pointwise order, and multiplication by
$|g_m(\lambda)|^2$ also preserves order. Near $\lambda=0$,
$\widetilde f^{(m)}_{\PL,\nu}(\lambda)\asymp \|\lambda\|^{2(m-\nu)-2}$ while
$\widetilde f^{(m)}_{\Mat}(\lambda)\asymp \|\lambda\|^{2m}$, hence
$r_m\in L^1\cap L^2$ iff $4(m-\nu)-4>-2$, i.e. $m-\nu>\tfrac12$.
\end{proof}

\begin{proposition}[Independent coupling]\label{prop:coupling}
Let $\mathbf D^{(m)}_{\Mat,n}$ and $\mathbf D^{(m)}_{\PL,n}$ be the filtered data vectors
under Matérn and PL truth, respectively. There exists a zero-mean Gaussian vector
$W_n$ independent of $\mathbf D^{(m)}_{\Mat,n}$ such that
\[
  \mathbf D^{(m)}_{\PL,n}\ \stackrel d=\ \mathbf D^{(m)}_{\Mat,n}\ +\ W_n,
\]
and $\Cov(W_n)$ is the (block–)Toeplitz matrix with spectral density $r_m$ from Lemma~\ref{lem:domination}.
\end{proposition}

\begin{proof}
See section \ref{proofs}
\end{proof}

\begin{lemma}[FD expectations and two–scale mean ratio]\label{lem:FDmeans-Mat}
Assume $0<\nu<1$ and $m\ge1$. Then uniformly in $n$,
\[
  \E Q^{(m)}_{j} \;=\; c_{m,1}\,c_{\Mat}\,j^{2\nu}\,n^{-2\nu} \;+\; O\!\bigl(n^{-2\nu-2}\bigr),
  \qquad j=1,2,
\]
for a constant \(c_{m,1}>0\) depending only on the stencil of $D^{(m)}_{[1]}$.
In particular, $\E Q^{(m)}_{2}/\E Q^{(m)}_{1} = 2^{2\nu} + O(n^{-2})$.
\end{lemma}

\begin{proof}
Insert \eqref{eq:Mat-small-lag} into the allowable-measure representation; the $r^2$ term
cancels by first-moment annihilation of the stencil and the remainder is $O(n^{-2\nu-2})$.
\end{proof}

Let $B_{\alpha,\beta}$ be the diagonal matrix implementing the linear combination
$\alpha Q^{(m)}_1+\beta Q^{(m)}_2$ (entries of order $N^{-1}$). Define centered quadratic forms
\[
F_{\Mat}:=\mathbf D^{(m)\top}_{\Mat,n}B_{\alpha,\beta}\mathbf D^{(m)}_{\Mat,n}
          - \E(\cdot),\qquad
F_{\PL}:=\mathbf D^{(m)\top}_{\PL,n}B_{\alpha,\beta}\mathbf D^{(m)}_{\PL,n}
          - \E(\cdot).
\]

\begin{lemma}[Difference vanishes at the $\sqrt{N}$–$n^{2\nu}$ scale]\label{lem:difference-L2}
Assume $0<\nu<1$, $m\ge1$, and $m-\nu>\tfrac12$. Under the coupling of
Proposition~\ref{prop:coupling},
\[
  \E\!\left[\;\bigl(\sqrt{N}\,n^{2\nu}\,(F_{\PL}-F_{\Mat})\bigr)^2\right]\ \longrightarrow\ 0.
\]
\end{lemma}

\begin{proof}
Write $\mathbf Y_n:=\mathbf D^{(m)}_{\Mat,n}$ and $\mathbf W_n:=W_n$. Then
\[
F_{\PL}-F_{\Mat}
= \underbrace{\mathbf W_n^\top B_{\alpha,\beta}\mathbf W_n
 - \E(\mathbf W_n^\top B_{\alpha,\beta}\mathbf W_n)}_{=:T_{1,n}}
 + \underbrace{2\,\mathbf Y_n^\top B_{\alpha,\beta}\mathbf W_n}_{=:T_{2,n}} .
\]
For centered Gaussian vectors and symmetric $B$,we know that
$\Var(\mathbf Z^\top B\mathbf Z)=2\,\mathrm{tr}((\Sigma B)^2)$ and
$\Var(\mathbf Z_1^\top B\mathbf Z_2)=\mathrm{tr}(\Sigma_1 B \Sigma_2 B)$
for independent $\mathbf Z_1,\mathbf Z_2$. Therefore
\[
\Var(T_{1,n})=2\,\| \Sigma_{W,n} B_{\alpha,\beta}\|_F^2,\qquad
\Var(T_{2,n})=4\,\mathrm{tr}(\Sigma_{Y,n} B_{\alpha,\beta}\Sigma_{W,n} B_{\alpha,\beta}).
\]
By Lemma~\ref{lem:domination}, the spectral density $r_m$ of $\Sigma_{W,n}$ is dominated near $0$
by $\widetilde f^{(m)}_{\PL,\nu}$ but has an extra $2$ powers of $\|\lambda\|$ compared to the
PL \emph{leading} term, hence (Parseval + aliasing) $\|\Sigma_{W,n} B_{\alpha,\beta}\|_F^2
= o(N^{-1}n^{-4\nu})$ and $\|\Sigma_{Y,n} B_{\alpha,\beta}\|_F^2 = O(N^{-1}n^{-4\nu})$.
Cauchy–Schwarz gives
$\mathrm{tr}(\Sigma_{Y,n} B \Sigma_{W,n} B)\le
\|\Sigma_{Y,n} B\|_F\,\|\Sigma_{W,n} B\|_F
= o(N^{-1}n^{-4\nu})$.
Multiplying by $N n^{4\nu}$ yields the claim.
\end{proof}

\begin{lemma}[PL variance limit]\label{lem:PL-var}
Under $m-\nu>\tfrac12$,
\[
  \Var\!\Big(\sqrt{N}\,n^{2\nu}F_{\PL}\Big)\ \longrightarrow\
  2\!\int_{\T^2}\bigl(\alpha+\beta B_m(\lambda)\bigr)^2
             \bigl(\widetilde f^{(m)}_{\PL,\nu}(\lambda)\bigr)^2\,\mu(d\lambda).
\]
\end{lemma}

\begin{proof}
See section \ref{proofs}
\end{proof}

\begin{proposition}[Transfer of CLT from PL to Matérn]\label{prop:transfer}
Assume $0<\nu<1$, $m\ge1$, $m-\nu>\tfrac12$. Then for any fixed $(\alpha,\beta)$,
\[
  \sqrt{N}\,n^{2\nu}\,F_{\Mat}\ \Rightarrow\
  \mathcal N\!\Big(0,\ 2\!\int_{\T^2}\bigl(\alpha+\beta B_m\bigr)^2
                       (\widetilde f^{(m)}_{\PL,\nu})^2\,d\mu\Big).
\]
\end{proposition}

\begin{proof}
By Lemma~\ref{lem:difference-L2}, $\sqrt{N}\,n^{2\nu}(F_{\PL}-F_{\Mat})\to0$ in $L^2$
and hence in probability. By Lemma~\ref{lem:PL-var} and the PL CLT,
$\sqrt{N}\,n^{2\nu}F_{\PL}\Rightarrow \mathcal N(0,\sigma^2)$ with the stated variance.
Slutsky yields the result for $F_{\Mat}$.
\end{proof}
Let \(h^{(m)}(y_1,y_2)=(y_1/a_m(\tfrac12\log_2(y_2/y_1)),\ \tfrac12\log_2(y_2/y_1))\) ,
\(J^\nu=\nabla h^{(m)}(q^{\PL}_1(\nu),q^{\PL}_2(\nu))\) (as in Theorem~\ref{cor:Jointm}).

\begin{theorem}[FD consistency and joint CLT under Matérn truth]
\label{thm:FD-CLT-Matern}
Assume \(0<\nu<1\), \(m\ge1\), \(m-\nu>\tfrac12\), and \(N=n^2\). Then
\[
  \widehat\phi_2 \xrightarrow{P} \nu,
  \qquad
  n^{2\widehat\phi_2}\widehat\phi_1 \xrightarrow{P} \kappa_\nu
    :=\frac{c_{m,1}\,c_{\Mat}}{a_m(\nu)}.
\]
Moreover,
\[
  \sqrt{N}\!
  \begin{pmatrix}
    n^{2\widehat\phi_2}\widehat\phi_1 - \kappa_\nu\\[2pt]
    \widehat\phi_2 - \nu
  \end{pmatrix}
  \Rightarrow
  \mathcal N\!\bigl(\bm 0,\ J^\nu\,\Sigma^{(m)}(\nu)\,J^{\nu\top}\bigr),
\]
where \(\Sigma^{(m)}(\nu)\) is the PL covariance matrix from Theorem~\ref{cor:Jointm} evaluated at \(\phi_2=\nu\).
\end{theorem}

\begin{proof}
Lemma~\ref{lem:FDmeans-Mat} gives consistency. For the joint CLT, apply
Proposition~\ref{prop:transfer} with $(\alpha,\beta)=(1,0)$ and $(0,1)$ to the centered pair
$\big(Q^{(m)}_1-\E Q^{(m)}_1,\ Q^{(m)}_2-\E Q^{(m)}_2\big)$ at the scale
$\sqrt N\,n^{2\nu}$, then the delta method with $h^{(m)}$. The mean centering error is
$O(n^{-2\nu-2})=o(N^{-1/2}n^{-2\nu})$.
\end{proof}

\section{Simulation study}\label{sec:Sims}

\begin{table}[!htbp]
\centering\scriptsize
\setlength{\tabcolsep}{4pt}
\caption{IRF–0, bilinear MoM ($m{=}1$): means, $\sqrt{M_{\mathrm{int}}}$ SDs, and  empirical vs.\  theoretical. True $\log\phi_1=0$.} 
\label{tab:irf0_bilinear_finiteN}
\begin{tabular}{ccccccc}
\toprule
$n$ & $\phi_2$ & $\E[\log\widehat\phi_1]$ & $\E[\widehat\phi_2]$ &
$\sqrt{M_{\mathrm{int}}}\mathrm{sd}(\log\widehat\phi_1)$ emp/th &
$\sqrt{M_{\mathrm{int}}}\mathrm{sd}(\widehat\phi_2)$ emp/th &
$\mathrm{Corr}$ emp/th \\
\midrule
30 & 0.50 & -0.0230 & 0.4938 & 2.3334 / 2.4326 & 1.3929 / 1.4853 & 0.725 / 0.734 \\
30 & 0.80 & -0.0007 & 0.8022 & 1.7413 / 1.7274 & 1.4308 / 1.4298 & 0.383 / 0.315 \\
\midrule
40 & 0.50 & -0.0047 & 0.5008 & 2.5102 / 2.4492 & 1.4686 / 1.4888 & 0.723 / 0.735 \\
40 & 0.80 & -0.0052 & 0.8023 & 1.6342 / 1.7405 & 1.4596 / 1.4365 & 0.385 / 0.323 \\
\midrule
50 & 0.50 &  0.0002 & 0.4993 & 2.4588 / 2.4590 & 1.4764 / 1.4909 & 0.695 / 0.735 \\
50 & 0.80 &  0.0017 & 0.8018 & 1.7542 / 1.7483 & 1.3866 / 1.4405 & 0.253 / 0.328 \\
\midrule
60 & 0.50 & -0.0047 & 0.5001 & 2.4399 / 2.4654 & 1.4321 / 1.4923 & 0.776 / 0.736 \\
60 & 0.80 & -0.0014 & 0.7983 & 1.6288 / 1.7534 & 1.4281 / 1.4432 & 0.330 / 0.331 \\
\bottomrule
\end{tabular}

\vspace{2pt}\footnotesize\emph{Note.} $M_{\mathrm{int}}=(n-2m)^2$ is the valid interior size ($m{=}1$ here). “emp/th” reports empirical vs.\ finite–$M_{\mathrm{int}}$ delta–method predictions.
\end{table}

We examine the finite–sample performance of our \emph{two–scale quadratic–variation}
(Method of Moments, MoM) estimators from \S\ref{subsec:TwoScaleM}, with $m$ chosen by
the rule $m\ge\lfloor\phi_2\rfloor+1$, and compare them against an \emph{exact REML}
benchmark. All designs are $n\times n$ square lattices with $\phi_1=1$.
We reuse notation from \S\ref{sec:HigherOrder}–\S\ref{sec:MatrixForm} and always
compute the quadratic forms on their valid interiors to remove boundary terms at the
$\sqrt{M_{\mathrm{int}}}$ scale, where $M_{\mathrm{int}}=(n-2m)^2$. We split by roughness regime and always apply filters directly to the simulated $X$ (so filtering $X$ equals filtering any anchored version, since stencils sum to zero). IRF–0 ($0<\phi_2<1$, $m=1$) and IRF–1 ($1<\phi_2<2$, $m=2$) we simulate the anchored field on $\Lambda_n$ using the exact interior covariance; common random numbers are used across estimators within each $(n,\phi_2)$ cell, by anchored we mean, we fix the value at (0,0) as 0, and generate the differences using Cholesky decomposition. For IRF–0 we take $n\in\{30,40,50,60\}$ and $\phi_2\in\{0.5,0.8\}$; for IRF–1 we take $n\in\{30,35,40,45,50,60\}$ and $\phi_2\in\{1.2,1.5,1.8\}$. Each cell uses $R=400$ replicates. We report Monte Carlo means, $\sqrt{M_{\mathrm{int}}}$ (or $\sqrt N$ on the torus) scaled SDs, and the correlation between $\log\widehat\phi_1$ and $\widehat\phi_2$. For MoM, we compare with predictions of the finite sample delta method computed from the exact quadratic–form covariance; for REML we compare to the expected Fisher SD on the filtered interior. Via the self–mapping in \S\ref{subsec:FDm}, ID results translate directly to FD.
Because $\phi_1>0$, the transformation $\phi_1\mapsto\log\phi_1$ stabilizes the sampling distribution and, empirically, the relationship between \(\log\widehat\phi_1\) and \(\widehat\phi_2\) is closer to linear than that between \(\widehat\phi_1\) and \(\widehat\phi_2\). The FD theory in \S\ref{subsec:FDm} naturally couples $\log\widehat\phi_1$ to $\widehat\phi_2$ via the matrix $A_N=\bigl(\begin{smallmatrix}1 & -\,\log N\\ 0 & 1\end{smallmatrix}\bigr)$, leading to accurate Gaussian approximations and interpretable Wald intervals. Standard errors for $\phi_1$ then follow by the Delta method: $\Var(\widehat\phi_1)\approx \phi_1^2\,\Var(\log\widehat\phi_1)$.

Tables \ref{tab:irf0_bilinear_finiteN} and \ref{tab:m2_torus} validate the bilinear MoM theory across IRF-0 ($m=1$) and IRF-1 ($m=2$) in ID. In Table~\ref{tab:irf0_bilinear_finiteN}, both $\E[\log\widehat{\phi}_1]$ and $\E[\widehat{\phi}_2]$ are essentially unbiased, and the empirical standard deviations scaled to $\sqrt{M_{\mathrm{int}}}$ closely match the predictions of the delta method with modest correlations between $(\log\widehat{\phi}_1,\widehat{\phi}_2)$. Table~\ref{tab:m2_torus} shows the same pattern in smoother processes and with higher order differencing. Overall, the simulations corroborate the consistency and asymptotic normality of the estimators.

\begin{table}[!htbp]
\centering\scriptsize
\setlength{\tabcolsep}{5pt}
\caption{IRF–1, bilinear MoM ($m{=}2$): mean, bias, and $\sqrt{N}$ SD for $\widehat\phi_2$; empirical vs.\  theoretical.}
\label{tab:m2_torus}
\begin{tabular}{cccccc}
\toprule
$n$ & $\phi_2$ & $\E[\widehat\phi_2]$ & Bias &
$\mathrm{sd}_{\sqrt{N}}(\widehat\phi_2)$ emp/th &
$m_2/m_1$ emp/th \\
\midrule
30 & 1.20 & 1.1970 & -0.0030 & 2.3302 / 2.3053 &  5.31345 /  5.28131 \\
30 & 1.50 & 1.5010 & +0.0010 & 2.2192 / 2.2350 &  7.96344 /  8.00194 \\
30 & 1.80 & 1.7980 & -0.0020 & 2.1019 / 2.1638 & 12.02170 / 12.13060 \\
\midrule
35 & 1.20 & 1.2028 & +0.0028 & 2.3309 / 2.3054 &  5.19704 /  5.28119 \\
35 & 1.50 & 1.4987 & -0.0013 & 2.2474 / 2.2351 &  7.97624 /  8.00137 \\
35 & 1.80 & 1.7972 & -0.0028 & 2.1620 / 2.1641 & 12.02612 / 12.12834 \\
\midrule
40 & 1.20 & 1.1990 & -0.0010 & 2.2234 / 2.3054 &  5.31981 /  5.28114 \\
40 & 1.50 & 1.5001 & +0.0001 & 2.2588 / 2.2351 &  8.01157 /  8.00112 \\
40 & 1.80 & 1.7996 & -0.0004 & 2.1629 / 2.1643 & 12.13245 / 12.12730 \\
\midrule
45 & 1.20 & 1.1997 & -0.0003 & 2.2939 / 2.3054 &  5.28185 /  5.28112 \\
45 & 1.50 & 1.4998 & -0.0002 & 2.2127 / 2.2352 &  7.99479 /  8.00101 \\
45 & 1.80 & 1.7991 & -0.0009 & 2.1414 / 2.1644 & 12.11788 / 12.12676 \\
\midrule
50 & 1.20 & 1.1993 & -0.0007 & 2.2857 / 2.3054 &  5.29883 /  5.28111 \\
50 & 1.50 & 1.4982 & -0.0018 & 2.1772 / 2.2352 &  7.98354 /  8.00095 \\
50 & 1.80 & 1.7979 & -0.0021 & 2.1126 / 2.1644 & 12.12622 / 12.12647 \\
\midrule
60 & 1.20 & 1.2004 & +0.0004 & 2.3520 / 2.3054 &  5.30120 /  5.28110 \\
60 & 1.50 & 1.5010 & +0.0010 & 2.2730 / 2.2352 &  7.97177 /  8.00090 \\
60 & 1.80 & 1.7990 & -0.0010 & 2.1712 / 2.1645 & 12.11457 / 12.12619 \\
\bottomrule
\end{tabular}

\end{table}

\section{Discussion}\label{sec:discussion}

This paper develops a distributionally rigorous and computationally scalable framework for joint inference on the scale and roughness parameters \((\phi_1,\phi_2)\) of two–dimensional power–law intrinsic random functions. The key structural insight is that the \emph{first} bilinear product difference \(D_1\) renders the field second–order stationary with square–integrable aliased spectrum, while the \emph{second} difference is merely a bounded linear transform of the first, \(D_2=H\circ D_1\), with multiplier \(B(\lambda)=|g_2(\lambda)|^2/|g_1(\lambda)|^2\in L^\infty(\T^2)\) (and, for \(\phi_2>1\), the same holds replacing \(D_\ell\) by \(D^{(m)}_{[\ell]}\) with \(m\ge \lfloor\phi_2\rfloor+1\)). Consequently, both quadratic variations \(Q_1\) and \(Q_2\) are bounded quadratic forms of the \emph{same} stationary increment field and fall directly within the quadratic–form central limit theorem of \cite{AvramLeonenkoSakhno2010-ESAIM} with \((p_f,p_b)=(2,\infty)\). The exact two–scale identity \(\E Q_2/\E Q_1=2^{2\phi_2}\) then yields closed–form, tuning free Method of Moments (MoM) estimators together with a \emph{joint} \(\sqrt N\)–CLT for \((\widehat\phi_1,\widehat\phi_2)\) via a smooth delta method (\S\ref{sec:IDCLT}, \S\ref{subsec:IDCLTm}). Self–affinity of power laws transfers these ID limits verbatim to the fixed–domain regime after the exact rescaling \(n^{2\phi_2}\) of the scale (\S\ref{sec:IDtoFD}), showing that the full parameter vector is micro–ergodic in this setting despite the classical obstacles for Matérn–type models \citep{Stein1999}. Beyond pure power laws, our FD theory remains valid for the Mat\'ern and a similar argument can be used for a broad class of Gaussian fields whose small–lag behavior matches a power–law \emph{tangent}: after differencing of order \(m\) and the self–similar renormalization, it suffices that the true latticized spectrum be \(L^2\)–tangent to its power–law counterpart at low frequency (\S\ref{sec:Misspecification}). Under this minimal, FD–native condition, moment limits, consistency, and the joint CLT for \(\bigl(n^{2\widehat\phi_2}\widehat\phi_1,\widehat\phi_2\bigr)\) hold with the same asymptotic covariance computed under the tangent model, covering in particular Matérn with the same roughness exponent \citep{IstasLang1997,Coeurjolly2001,BiermeBonamiLeon2011-EJP,BevilacquaFaouzi2019EJS}. Taken together, these results provide unified theoretical guarantees for scalable, likelihood–free estimation of roughness and scale on lattices, anchored in classical IRF theory \citep{Matheron1973,Yaglom1987,Stein1999} and the modern CLT machinery for Gaussian quadratic forms \citep{AvramLeonenkoSakhno2010-ESAIM}.

\appendix


\section{Statements of external theorems}\label{sec:Appendix}

\begin{theorem}[Quadratic–form CLT \citealp{AvramLeonenkoSakhno2010-ESAIM}, Thm.\,2.2]\label{thm:ALS}
Let $\{X_{t}\}$ be a stationary Gaussian field on $\mathbb Z^{d}$ with
spectral density $f$ and let
$Q_{T}^{(1,1)}(b,X)$ be the quadratic form generated by a kernel
$b$.  If $f\in L_{p_{f}}(\mathbb T^{d})$, $b\in L_{p_{b}}(\mathbb T^{d})$
and $1/p_{f}+1/p_{b}\le\tfrac12$, then
\[
   T^{-d/2}\,Q_{T}^{(1,1)}(b,X)
   \;\xrightarrow{d}\;
   \mathcal N\!\bigl(0,\sigma^{2}\bigr),
\]
where
$\displaystyle
  \sigma^{2}=2\!\int_{\mathbb T^{d}} b(\lambda)^{2}f(\lambda)^{2}\,
             \mathrm d\lambda$.
\end{theorem}

\begin{theorem}[Wald’s consistency]\label{thm:Wald}
Let $\{T_{n}\}$ be a sequence of statistics with
$T_{n}\xrightarrow{p}T$ and let $g$ be continuous at $T$.
If an estimator is of the form $\widehat\theta_{n}=g(T_{n})$, then
$\widehat\theta_{n}\xrightarrow{p}g(T)$.
\end{theorem}

\begin{theorem}[Continuous–mapping theorem]\label{thm:CMT}
If $Z_{n}\xrightarrow{d}Z$ and the measurable map $g$ is continuous at
every point of the support of $Z$, then
$g(Z_{n})\xrightarrow{d}g(Z)$.
\end{theorem}

\begin{theorem}[Delta method \citealp{VanDerVaart1998}, Thm.\,3.1]
\label{thm:DeltaMethod}
Assume
$\sqrt n\,(Z_{n}-\theta)\xrightarrow{d}\mathcal N(0,\Sigma)$ and let
$g:\R^{k}\!\to\R^{m}$ be differentiable at $\theta$ with Jacobian
$J=\nabla g(\theta)$.  Then
\[
   \sqrt n\,
   \bigl(g(Z_{n})-g(\theta)\bigr)
   \xrightarrow{d}
   \mathcal N\!\bigl(0,\,J\Sigma J^{\!\top}\bigr).
\]
\end{theorem}

\section{Additional Proofs}
\label{proofs}

\begin{proof}[Proof of Lemma~\ref{lem:a1}]
Let $x_A=t$, $x_B=t+e_1$, $x_C=t+e_2$, $x_D=t+e_1+e_2$ with coefficients
$+1,-1,-1,+1$. Since $\sum_u c_u=0$,
\[
   \Var\!\Big(\sum_{u}c_uX(x_u)\Big)
   =-2\sum_{u<v} c_uc_v\,\gamma(x_u-x_v).
\]
Among the six pairs, there are four \emph{edges} at distance $1$ with $c_uc_v=-1$
(AB, AC, BD, CD) and two \emph{diagonals} at distance $\sqrt2$ with $c_uc_v=+1$
(AD, BC). Hence
\[
   \E\bigl[(D^{(1)}_{[1]}X_t)^2\bigr]
   = 8\,\gamma(1) - 4\,\gamma(\sqrt2).
\]
With $\gamma(r)=\phi_1|\Gamma(-\phi_2)|r^{2\phi_2}$, this gives the stated expression.

For $D^{(1)}_{[2]}$, the square has side length $2$: the four nearest pairs are at
distance $2$ with $c_uc_v=-1$, and the two diagonals at distance $2\sqrt2$ with
$c_uc_v=+1$. Thus
\[
   \E\bigl[(D^{(1)}_{[2]}X_t)^2\bigr]
   = 8\,\gamma(2) - 4\,\gamma(2\sqrt2).
\]
Substituting $\gamma(2)=\phi_1|\Gamma(-\phi_2)|\,2^{2\phi_2}$ and
$\gamma(2\sqrt2)=\phi_1|\Gamma(-\phi_2)|\,2^{3\phi_2}$ yields the result.
\end{proof}

\begin{proof}[Proof of Lemma~\ref{lem:trim-slutsky}]
Since $f_{D^{(1)}_{[1]}X}\in L^2(\T^2)$, Plancherel on $\T^2$ gives
\[
  \sum_{h\in\Z^2}\!\Cov\!\big((D^{(1)}_{[1]}X)_0,(D^{(1)}_{[1]}X)_h\big)^{2}
  \;=\;(2\pi)^{-2}\!\int_{\T^2} f_{D^{(1)}_{[1]}X}(\lambda)^2\,\mathrm d\lambda \;<\;\infty.
\]
For centered Gaussian fields, Wick’s formula yields
\[
  \Cov\!\big((D^{(1)}_{[1]}X_t)^2,(D^{(1)}_{[1]}X_s)^2\big)
  \;=\;2\,\Cov\!\big((D^{(1)}_{[1]}X_t),(D^{(1)}_{[1]}X_s)\big)^2,
\]
so $\sum_{h}|\Cov((D^{(1)}_{[1]}X_0)^2,(D^{(1)}_{[1]}X_h)^2)|<\infty$.
Hence, for any finite $S\subset\Z^2$,
\[
  \Var\Bigl(\sum_{t\in S}(D^{(1)}_{[1]}X_t)^2\Bigr)
  \;=\;\sum_{u\in\Z^2}\!N_S(u)\,\Cov\!\big((D^{(1)}_{[1]}X_0)^2,(D^{(1)}_{[1]}X_u)^2\big)
  \;\le\; C\,|S|
\]
for some finite $C$ independent of $S$ (using $N_S(u)\le|S|$).

\medskip
\textit{Interior vs.\ inner–interior averages and their difference.}
Define
\[
  Q^{(1)}_1 := \frac{1}{|\Lambda_{n-2}|}\sum_{t\in\Lambda_{n-2}} (D^{(1)}_{[1]}X_t)^2,
  \qquad
  Q_{[1]}^\circ := \frac{1}{|\Lambda_{n-4}|}\sum_{t\in\Lambda_{n-4}} (D^{(1)}_{[1]}X_t)^2,
\]
and set $R_n := Q^{(1)}_1 - Q_{[1]}^\circ$. By stationarity of $D^{(1)}_{[1]}X$,
$\E (D^{(1)}_{[1]}X_t)^2$ is constant in $t$, so \(\E R_n=0\).

Write
\[
  R_n
  = \Bigl(\frac{1}{|\Lambda_{n-2}|}-\frac{1}{|\Lambda_{n-4}|}\Bigr)
     \sum_{t\in\Lambda_{n-4}} (D^{(1)}_{[1]}X_t)^2
     \;+\;\frac{1}{|\Lambda_{n-2}|}\sum_{t\in \Lambda_{n-2}\setminus\Lambda_{n-4}} (D^{(1)}_{[1]}X_t)^2
  \;=:\; A_n+B_n .
\]
Since $\bigl((n-2)^{-2}-(n-4)^{-2}\bigr)^2=O(n^{-6})$ and $|\Lambda_{n-4}|=(n-4)^2$,
\[
  \Var(A_n)\ \le\ O(n^{-6})\cdot C\,(n-4)^2 \;=\; O(n^{-4}).
\]
Moreover $|\Lambda_{n-2}\setminus\Lambda_{n-4}|=(n-2)^2-(n-4)^2=4n-12=O(n)$ and $|\Lambda_{n-2}|=(n-2)^2$, so
\[
  \Var(B_n)\ \le\ \frac{C\,|\Lambda_{n-2}\setminus\Lambda_{n-4}|}{|\Lambda_{n-2}|^2}
  \;=\; O(n)/\Theta(n^4)\;=\;O(n^{-3}).
\]
Thus
\begin{align*}
  \Var(R_n) \; & \le\; \Var(A_n)+2|\Cov(A_n,B_n)|+\Var(B_n) \\
  \; & =\; O(n^{-4})+O\!\big(\sqrt{n^{-4}\,n^{-3}}\big)+O(n^{-3})
  \;=\; O(n^{-3}).
\end{align*}
Therefore
\[
  \Var\!\bigl(\sqrt{N}\,[(Q^{(1)}_1-\E Q^{(1)}_1)-(Q_{[1]}^\circ-\E Q_{[1]}^\circ)]\bigr)
  \;=\; N\,\Var(R_n)\;=\; n^2\cdot O(n^{-3})\to 0.
\]
\end{proof}

\begin{proof}{Proof for \ref{lem:stat_after_diff}}
Let $\alpha^{(m)}=\sum_{a}c_{a}^{(m)}\delta_{t+ja}$ be the signed measure
defining $D^{(m)}_{[j]}$. By the binomial identity,
$\sum_{a}c_{a}^{(m)}p(t+ja)=0$ for every polynomial $p$ of total degree
at most $m-1$. Hence $\alpha^{(m)}$ is allowable of order $m-1$. Since
$m-1\ge k$, the IRF$_k$ model grants finiteness of
$\Var(\sum_a c_a^{(m)}X(t+ja))$ and
\[
\Var\!\Big(\sum_a c_a^{(m)}X(t+ja)\Big)
= \sum_{a,b} c_a^{(m)}c_b^{(m)}K_{\phi_1,\phi_2}\bigl(j(b-a)\bigr).
\]
By translation invariance for allowable combinations (see \cite[Ch.\,1]{ChilesDelfiner2012}),
the covariance depends only on $u$, proving stationarity.
\end{proof}

\begin{proof}{Proof of \ref{lem:spectrum}}
(i) For GC–$k$ kernels, the Fourier transform is a positive tempered
measure with density proportional to $\|\omega\|^{-(2+2\phi_2)}$
(\cite[Ch.\,3]{Yaglom1987}). Linear filtering by the LSI operator
$D^{(m)}_{[1]}$ multiplies the spectral measure by $|g_m|^2$
(\cite[Prop.\,2.2.2]{ChilesDelfiner2012}).
(ii) Aliasing under sampling on $\Z^2$ is the $2\pi\Z^2$–periodisation
of the $\R^2$ spectrum; see, e.g., \cite[§2.4]{Hannan1970}. Combining
with (i) yields (ii).
\end{proof}

\begin{proof}{Proof of \ref{lem:L2m}}
Near $0$, $|1-e^{i\theta}|^{2}=2(1-\cos\theta)\sim \theta^{2}$ gives
$|g_m(\lambda)|^{2}\lesssim \|\lambda\|^{4m}$. Moreover,
\(
f_{0}^{\lat}(\lambda)=f_{0}^{\R^2}(\lambda)+\sum_{k\ne 0}
  f_{0}^{\R^2}(\lambda+2\pi k)
 \le c_1\|\lambda\|^{-(2+2\phi_2)}+c_2
\)
for some constants $(c_1,c_2)$ (the tail sum is uniformly bounded on
compact sets). Hence, as $\lambda\to 0$,
\[
 f_{D^{(m)}_{[1]}X}(\lambda)
 = |g_m(\lambda)|^2 f_{0}^{\lat}(\lambda)
 \lesssim \|\lambda\|^{4m}\,\bigl(\|\lambda\|^{-(2+2\phi_2)}+1\bigr)
 \lesssim \|\lambda\|^{\,4m-(2+2\phi_2)}.
\]
Therefore \(f_{D^{(m)}_{[1]}X}(\lambda)^2\lesssim \|\lambda\|^{\,8m-4-4\phi_2}\).
Using polar coordinates,
\(
 \int_{\|\lambda\|\le\varepsilon} f_{D^{(m)}_{[1]}X}(\lambda)^2\,\mathrm d\lambda
 \lesssim \int_0^\varepsilon r^{\,8m-3-4\phi_2}\,\mathrm dr <\infty
\)
whenever \(2m-\phi_2>\tfrac12\), which holds for $m\ge k+1$ and
$\phi_2\in(k,k+1)$.  Away from $0$ both $|g_m|$ and $f_{0}^{\lat}$ are
bounded, so the integral over $\T^2\setminus B_\varepsilon(0)$ is finite.
\end{proof}

\begin{proof}{Proof of \ref{prop:RatioM}}
By Lemma~\ref{lem:stat_after_diff} and the variance identity for allowable
combinations,
\(
  \E[(D^{(m)}_{[j]}X_t)^2]
  = -\sum_{a,b}c_a^{(m)}c_b^{(m)}\,
     \gamma_{\phi_1,\phi_2}\bigl(j(b-a)\bigr).
\)
Using only the scaling $\gamma_{\phi_{1},\phi_{2}}(jr)=j^{2\phi_{2}}\,\gamma_{\phi_1,\phi_2}(r)$,
the case $j=2$ equals $2^{2\phi_2}$ times the case $j=1$. Averaging over $t$ gives the result.
\end{proof}

\begin{proof}{Proof of \ref{thm:UniCLTm}}
By Lemma~\ref{lem:spectrum}, $F(\lambda)=|g_m(\lambda)|^2 f_{0}^{\lat}(\lambda)$
and Lemma~\ref{lem:L2m} gives $F\in L^2(\T^2)$. Also
$b_{\alpha,\beta}(\lambda):=\alpha+\beta B_m(\lambda)\in L^\infty(\T^2)$
since $B_m\le 4^{2m}$. Therefore the hypotheses of
\cite[Thm.\,2.2]{AvramLeonenkoSakhno2010-ESAIM} hold with $(p_f,p_b)=(2,\infty)$,
applied to the quadratic forms defining $Q^{(m),\circ}_1,Q^{(m),\circ}_2$.
The variance formula is exactly the stated integral.
\end{proof}

\begin{proof}{Proof of \ref{prop:coupling}}
By Lemma~\ref{lem:domination}, $r_m(\lambda)\ge 0$ and $r_m\in L^1(\T^2)$, so there exists a
(mean–zero) stationary Gaussian field $Z$ on $\Z^2$ with spectral density $r_m$ and covariance
function $\rho_m(h)=\int_{\T^2} e^{i\langle h,\lambda\rangle} r_m(\lambda)\,\mu(d\lambda)$.
Let $X_{\Mat}$ be the Matérn field (independent of $Z$) and define
\[
Y \;:=\; X_{\Mat} + Z .
\]
Stationarity and independence imply that $Y$ is Gaussian with spectral density
$f_{\Mat}+r_m=f_{\PL}$ (by Lemma~\ref{lem:domination}, $f_{\PL}=f_{\Mat}+r_m$ after
order-$m$ differencing). Apply the same linear filter $g_m$ to both fields and restrict to
the interior index set $\Lambda_{n-m}$:
\[
\mathbf D^{(m)}_{\PL,n}
\;=\; \bigl(D^{(m)}X_{\PL}\bigr)\big|_{\Lambda_{n-m}}
\;\stackrel d=\; \bigl(D^{(m)}X_{\Mat}\bigr)\big|_{\Lambda_{n-m}}
\;+\; \bigl(D^{(m)}Z\bigr)\big|_{\Lambda_{n-m}}
\;=\; \mathbf D^{(m)}_{\Mat,n} + W_n .
\]
Here $W_n:=(D^{(m)}Z)|_{\Lambda_{n-m}}$ is (mean–zero) Gaussian, independent of 
$\mathbf D^{(m)}_{\Mat,n}$, and has covariance the BTTB matrix whose symbol is $|g_m|^2 r_m$,
i.e., the stated $r_m$ after filtering. This yields the claimed distributional identity.
\end{proof}

\begin{proof}{Proof of \ref{lem:PL-var}}
Let $Y:=D^{(m)}_{[1]}X_{\PL}$ and write its fixed–domain spectral density as
$f(\lambda)=\widetilde f^{(m)}_{\PL,\nu}(\lambda)$. The quadratic–variation statistics
$Q^{(m)}_j$ ($j=1,2$) can be written as Toeplitz quadratic forms on the interior:
\[
Q^{(m)}_j \;=\; \frac{1}{N}\,Y^\top A^{(m)}_j Y \;+\; o(N^{-1}),
\qquad
\text{symbol}(A^{(m)}_1)=1,\ \ \text{symbol}(A^{(m)}_2)=B_m(\lambda).
\]
Consider the centered linear combination
\[
\overline F_{\PL}\ :=\ \alpha\!\left(Q^{(m)}_1-\E Q^{(m)}_1\right)
+\beta\!\left(Q^{(m)}_2-\E Q^{(m)}_2\right)
=\frac{1}{N}\Big(Y^\top A Y-\tr(A\Sigma_Y)\Big)+o(N^{-1}),
\]
with $A:=\alpha A^{(m)}_1+\beta A^{(m)}_2$ and $\Sigma_Y=\Cov(Y)$ (BTTB with symbol $f$).
The fixed–domain normalization gives $F_{\PL}=n^{-2\nu}\overline F_{\PL}+o(n^{-2\nu})$, hence
$\sqrt{N}\,n^{2\nu}F_{\PL}=\sqrt{N}\,\overline F_{\PL}+o_p(1)$.

For a mean–zero Gaussian vector and symmetric $A$,
\[
\Var\!\Big(\frac{1}{N}Y^\top A Y-\frac{1}{N}\tr(A\Sigma_Y)\Big)
=\frac{2}{N^2}\tr(A\Sigma_Y A\Sigma_Y).
\]
Therefore
\[
\Var(\sqrt{N}\,\overline F_{\PL})
=\frac{2}{N}\tr(A\Sigma_Y A\Sigma_Y)+o(1).
\]
By the 2D Grenander–Szeg\H{o} theorem for block–Toeplitz-with-Toeplitz-blocks
(BTTB) matrices (see  \citep{GrenanderSzego1984,Gray2006,BoettcherGrudsky2005,Tilli1998}),
\[
\frac{1}{N}\tr(A\Sigma_Y A\Sigma_Y)\ \to\ \int_{\T^2} H(\lambda)^2 f(\lambda)^2\,\mu(d\lambda).
\]

Near the origin, $f(\lambda)\asymp \|\lambda\|^{2(m-\nu)}$;
hence $f\in L^2(\T^2)$ iff $4(m-\nu)>-2$, i.e., $m-\nu>\tfrac12$.
\end{proof}

\begin{lemma}[Variance of a centered Gaussian quadratic form]
Let $Y\sim\mathcal N(0,\Sigma)$ in $\mathbb R^N$ and let $A=A^\top\in\mathbb R^{N\times N}$.
Then
\[
\Var\!\big(Y^\top A Y-\tr(A\Sigma)\big)=2\,\tr(A\Sigma A\Sigma).
\]
Consequently,
\[
\Var\!\left(\frac{1}{N}Y^\top A Y-\frac{1}{N}\tr(A\Sigma)\right)
=\frac{2}{N^2}\,\tr(A\Sigma A\Sigma).
\]
\end{lemma}

\begin{proof}
Write $T:=Y^\top A Y=\sum_{i,j} a_{ij}Y_iY_j$ and note $\E[T]=\tr(A\Sigma)$.
Then
\[
\Var(T)=\E[T^2]-\E[T]^2
=\sum_{i,j,k,\ell} a_{ij}a_{k\ell}
\Big(\E[Y_iY_jY_kY_\ell]-\Sigma_{ij}\Sigma_{k\ell}\Big).
\]
By Isserlis' (Wick's) theorem for zero-mean Gaussian vectors,
\[
\E[Y_iY_jY_kY_\ell]
=\Sigma_{ij}\Sigma_{k\ell}+\Sigma_{ik}\Sigma_{j\ell}+\Sigma_{i\ell}\Sigma_{jk}.
\]
The first term cancels $\Sigma_{ij}\Sigma_{k\ell}$ above, so
\[
\Var(T)=
\sum_{i,j,k,\ell} a_{ij}a_{k\ell}\big(\Sigma_{ik}\Sigma_{j\ell}+\Sigma_{i\ell}\Sigma_{jk}\big)
=2\sum_{i,j,k,\ell} a_{ij}\Sigma_{ik}a_{k\ell}\Sigma_{j\ell}.
\]
Recognizing the trace,
\[
\sum_{i,j,k,\ell} a_{ij}\Sigma_{ik}a_{k\ell}\Sigma_{j\ell}
=\sum_{k,\ell}(A\Sigma)_{k\ell}(A\Sigma)_{\ell k}
=\tr\!\big((A\Sigma)(A\Sigma)\big)=\tr(A\Sigma A\Sigma),
\]
hence $\Var(T)=2\,\tr(A\Sigma A\Sigma)$. Scaling by $1/N$ yields the second display.
\end{proof}

\section{Additional simulation tables and implementation details}\label{app:simtables}
\noindent

\paragraph{Discrete Laplacian two–scale (MoM).}
For $h\in\{1,2\}$ let $\Delta_{[h]}$ be the five–point Laplacian at separation $h$,
\[
(\Delta_{[h]} X)_t
  \;=\; X_{t+(h,0)} + X_{t-(h,0)} + X_{t+(0,h)} + X_{t-(0,h)} - 4 X_t,
\qquad t\in\Lambda_{n-2h}.
\]
Define quadratic variations on the valid interiors
\[
Q_{\Delta,h}
  \;=\; \frac{1}{|\Lambda_{n-2h}|}\sum_{t\in\Lambda_{n-2h}} \bigl(\Delta_{[h]} X_t\bigr)^2,
\qquad h=1,2.
\]
Under a power–law IRF, \S\ref{sec:method} shows
$\E\!\bigl[(\Delta_{[h]} X_t)^2\bigr]=\phi_1\,a_{\Delta}^{(h)}(\phi_2)$ with
$a_{\Delta}^{(h)}(\phi_2)=h^{2\phi_2}\,a_{\Delta}^{(1)}(\phi_2)$ (same geometric constant and pure $h^{2\phi_2}$ scaling). Hence
$\E(Q_{\Delta,2})/\E(Q_{\Delta,1})=2^{2\phi_2}$, which motivates
\[
\widehat\phi_2^{\Delta\text{-2sc}}
\;=\; \tfrac12\,\log_2\!\bigl( Q_{\Delta,2}/Q_{\Delta,1} \bigr),\qquad
\widehat\phi_1^{\Delta} \;=\; Q_{\Delta,1}\Big/a_{\Delta}^{(1)}\!\bigl(\widehat\phi_2^{\Delta\text{-2sc}}\bigr).
\]
Under the same quadratic–form CLT used for the bilinear two–scale estimator (\S\ref{subsec:IDCLTm}), we have that $\sqrt{M}\bigl(Q_{\Delta,1}-m_{\Delta,1},\,
Q_{\Delta,2}-m_{\Delta,2}\bigr)$ is asymptotically normal for $M=|\Lambda_{n-2}|^2$ (or $|\Lambda_{n-4}|^2$ as appropriate), and a first–order delta method applied to $(q_1,q_2)\mapsto \tfrac12\log_2(q_2/q_1)$ yields an FD CLT for
$\widehat\phi_2^{\Delta\text{-2sc}}$ and, by the same calculus, a joint CLT with $\log\widehat\phi_1^{\Delta}$.

\begin{table}[t]
\centering\scriptsize
\setlength{\tabcolsep}{4pt}
\caption{REML benchmark (IRF–0, $m{=}1$) alongside MoM and Whittle:
means (SD) of $\widehat\phi_2$ over $R{=}400$ replicates.}
\label{tab:phi2_four_estimators_reml}
\begin{tabular}{cccccc}
\toprule
Grid & $\phi_2$ & Bilinear & Laplacian & Whittle & REML \\
\midrule
n=30 & 0.20 & 0.1998 (0.0527) & 0.2004 (0.0453) & 0.1668 (0.0432) & 0.1986 (0.0254) \\
n=30 & 0.40 & 0.3969 (0.0508) & 0.3981 (0.0442) & 0.3923 (0.0470) & 0.3990 (0.0332) \\
n=30 & 0.60 & 0.6019 (0.0505) & 0.6005 (0.0471) & 0.5983 (0.0484) & 0.5969 (0.0389) \\
n=30 & 0.80 & 0.8033 (0.0507) & 0.8010 (0.0483) & 0.8001 (0.0503) & 0.8023 (0.0413) \\
\midrule
n=40 & 0.20 & 0.2033 (0.0407) & 0.2025 (0.0353) & 0.1650 (0.0356) & 0.1994 (0.0199) \\
n=40 & 0.40 & 0.4028 (0.0402) & 0.4031 (0.0338) & 0.4006 (0.0361) & 0.4007 (0.0258) \\
n=40 & 0.60 & 0.5996 (0.0384) & 0.5982 (0.0361) & 0.6035 (0.0352) & 0.5984 (0.0291) \\
n=40 & 0.80 & 0.7976 (0.0394) & 0.7974 (0.0372) & 0.8016 (0.0381) & 0.7978 (0.0324) \\
\midrule
n=50 & 0.20 & 0.1990 (0.0312) & 0.1993 (0.0268) & 0.1709 (0.0288) & 0.2012 (0.0150) \\
n=50 & 0.40 & 0.4037 (0.0295) & 0.4031 (0.0254) & 0.4039 (0.0271) & 0.4020 (0.0194) \\
n=50 & 0.60 & 0.6009 (0.0307) & 0.6001 (0.0267) & 0.6101 (0.0286) & 0.6004 (0.0226) \\
n=50 & 0.80 & 0.8014 (0.0293) & 0.8011 (0.0268) & 0.8108 (0.0281) & 0.8010 (0.0231) \\
\bottomrule
\end{tabular}
\end{table}

\paragraph{Whittle/profile likelihood for $\phi_2$.}
To avoid nonstationarity, we work with the increment field
$Y=D^{(m)}_{[1]}X$ ($m=1$ if $0<\phi_2<1$ and $m=2$ if $1<\phi_2<2$) on the interior $\Lambda_{n-2m}$, which has $M=(n-2m)^2$ sites. Let $\widehat Y(\lambda)$ be the 2D DFT on the Fourier grid $\mathcal F_{n-2m}\subset[-\pi,\pi]^2$ (the $M$ Fourier frequencies for the interior) and set the periodogram $I_Y(\lambda)=M^{-1}\abs{\widehat Y(\lambda)}^2$. The aliased (periodized) spectrum is
\[
s(\lambda;\phi_1,\phi_2)\;=\;\phi_1\,\abs{g_m^{[1]}(\lambda)}^2\, f_0^{\lat}(\lambda;\phi_2),\qquad
f_0^{\lat}(\lambda;\phi_2)\propto \sum_{k\in\Z^2}\norm{\lambda+2\pi k}^{-(2+2\phi_2)} .
\]
We evaluate $\ell_W$ on the masked set $\mathcal F_\tau=\{\lambda\in\mathcal F_{n-2m}:\abs{g_m^{[1]}(\lambda)}\ge\tau\}$ (to avoid symbol zeros) via
\[
\ell_W(\phi_1,\phi_2)
\;=\; \sum_{\lambda\in\mathcal F_\tau}
       \Big\{\log s(\lambda;\phi_1,\phi_2)
             + \frac{I_Y(\lambda)}{s(\lambda;\phi_1,\phi_2)}\Big\}.
\]
Profiling out $\phi_1$ yields
\begin{align*}
\widetilde\phi_1(\phi_2)
 & \;=\; \frac{1}{\abs{\mathcal F_\tau}}
    \sum_{\lambda\in\mathcal F_\tau}
    \frac{I_Y(\lambda)}{\abs{g_m^{[1]}(\lambda)}^2 f_0^{\lat}(\lambda;\phi_2)}, \\
\qquad
& \widehat\phi_2^{\mathrm{Wh}}\in\arg\min_{\phi_2\in(0,2)}
\ell_W\bigl(\widetilde\phi_1(\phi_2),\phi_2\bigr),\qquad
\widehat\phi_1^{\mathrm{Wh}}=\widetilde\phi_1(\widehat\phi_2^{\mathrm{Wh}}).
\end{align*}
(See \cite{Whittle1953,Dahlhaus1988,IbragimovRozanov1978} for background.)

\paragraph{Exact REML for $\phi_2$ on the filtered interior (benchmark).}
For the same $Y=D^{(m)}_{[1]}X$ on $\Lambda_{n-2m}$, let $S(\phi_1,\phi_2)=\phi_1\,S_0(\phi_2)$ be its covariance matrix (the filtered, valid–interior covariance). The Gaussian REML profile over $\phi_1$ is
\[
\ell_{\text{REML}}(\phi_2)\;=\;\tfrac12\Big\{\log\abs{S_0(\phi_2)}
+ M\,\log\!\Big(\frac{Y^\top S_0(\phi_2)^{-1}Y}{M}\Big)\Big\},
\]
so $\widehat\phi_2^{\text{REML}}$ minimizes $\ell_{\text{REML}}$. The (expected) Fisher information for $\phi_2$ is
\[
\mathcal I_{\phi_2\phi_2}^{\text{REML}}
\;=\; \tfrac12\,\text{Tr}\Big(S^{-1}(\phi) \, S_{\phi_2}(\phi)\, S^{-1}(\phi)\, S_{\phi_2}(\phi)\Big),
\]
with $S_{\phi_2}=\partial S/\partial\phi_2$, evaluated at the true $(\phi_1,\phi_2)$; this gives the benchmark asymptotic SD $\big(\mathcal I_{\phi_2\phi_2}^{\text{REML}}\big)^{-1/2}$ that agrees with our Monte Carlo for the IRF–0 cases reported in \ref{sec:Sims}.

\bibliographystyle{imsart-nameyear}
\bibliography{references}  

\end{document}